\documentclass[12pt]{article}
\usepackage{amssymb,amsmath}
\newcommand{\GL}{{\rm GL}}

\newcommand{\kb}{\overline{k}}
\renewcommand{\and}{\;\;\;\mbox{and}\;\;\;}
\newcommand{\modd}[1]{\; ( \text{mod} \; #1)}
\newcommand{\beql}[1]{\begin{equation}\label{#1}}
\newcommand{\eeq}{\end{equation}}
\renewcommand{\b}[1]{\mathbf{#1}}
\newcommand{\cl}[1]{{\cal #1}}
\newcommand{\cP}{\cl{P}}
\newcommand{\ov}{\mathcal{O}_v}
\newcommand{\ok}{\mathcal{O}_k}
\newcommand{\Z}{\mathbb{Z}}

\newcommand{\R}{\mathbb{R}}
\newcommand{\Q}{\mathbb{Q}}
\newcommand{\proj}{\mathbb{P}}
\newcommand{\x}{\b{x}}
\newcommand{\rank}{{\rm rank}}
\newcommand{\diag}{{\rm diag}}

\newcommand{\ep}{\varepsilon}
\newtheorem{prop}{Proposition}
\newtheorem{theorem}{Theorem}
\numberwithin{equation}{section}
\newtheorem{lemma}{Lemma}[section]
\newcommand{\hd}{{\rm det}_{/2}}
\begin{document}

\title{Zeros of Pairs of Quadratic Forms}
\author{D.R. Heath-Brown\\Mathematical Institute, Oxford}
\date{}
\maketitle

\section{Introduction}\label{sec-intr}

Let $Q_1(X_1,\dots,X_n)$ and $Q_2(X_1,\dots,X_n)$ be a pair of
quadratic forms defined over $\Q$, or more generally over a number
field $k$.  This paper will be concerned with the existence of 
non-trivial simultaneous zeros of the two forms, over $\Q$ or $k$ as
appropriate. In particular one may hope for a local-to-global
principle.  Indeed one may also ask whether the global points are
dense in the ad\`{e}lic points.  When $n\ge 9$ this problem was given
a very satisfactory
treatment in the work of Colliot-Th\'{e}l\`{e}ne, Sansuc and
Swinnerton-Dyer \cite{CTSSD1} and \cite{CTSSD2}.  
It follows from \cite[Theorem A(i)(a), page 40]{CTSSD1} that the forms have a
non-trivial common zero providing that they have a nonsingular common
zero in every completion of $k$, and from \cite[Theorem A(ii), page 40]{CTSSD1}
that the weak approximation principle holds if the intersection
$Q_1=Q_2=0$ is nonsingular.  Moreover they proved 
\cite[Theorem C, page 38]{CTSSD1} that when $n\ge 9$ there 
is automatically a non-trivial common zero if the field $k$ is 
totally imaginary.

Colliot-Th\'{e}l\`{e}ne, Sansuc and Swinnerton-Dyer \cite[Section 15]{CTSSD2}
give a number of examples showing that such results would be false for 
forms in sufficiently few variables.  In particular the example
\[Q_1=X_1X_2-(X_3^2-5X_4^2),\;\;\;
Q_2=(X_1+X_2)(X_1+2X_2)-(X_3^2-5X_5^2)\]
due to Birch and Swinnerton-Dyer \cite{BSD}, 
gives a smooth intersection violating the Hasse principle,
while the example
\[Q_1=X_1X_2-(X_3^2+X_4^2),\;\;\;
Q_2=(4X_2-3X_1)(4X_1-X_2)-(X_3^2+X_5^2)\]
due to Colliot-Th\'el\`ene and Sansuc \cite{CTS}, gives a smooth 
intersection for which weak approximation fails.
Colliot-Th\'{e}l\`{e}ne, Sansuc and Swinnerton-Dyer 
conjecture \cite[Section 16]{CTSSD2}
that the Hasse principle should hold for nonsingular
intersections as soon as $n\ge 6$. The goal of the present
paper is the corresponding result for forms in 8 variables.

\begin{theorem}\label{T1}
Let $Q_1(X_1,\dots,X_8)$ and 
$Q_2(X_1,\dots,X_8)$ be two quadratic forms over a number field $k$ 
such that the
projective variety 
\[\cl{V}:\; Q_1(\b{X})=Q_2(\b{X})=0\]
is nonsingular.  Then
the Hasse principle and weak approximation hold for $\cl{V}$.
\end{theorem}

For pairs of diagonal forms this was shown by Colliot-Th\'el\`ene, in
unpublished work.  The proof used the corresponding special case of
Theorem 2, which he established using the argument ascribed to
him in the proof of Lemma \ref{thelemma}.

Perhaps it is as well to explain precisely what we mean by the
variety $Q_1(X_1,\dots,X_n)=Q_2(X_1,\dots,X_n)=0$ being nonsingular
over a general field $K$.  By this it is meant that the matrix
\[\left(\begin{array}{ccc}
\frac{\partial Q_1}{\partial x_1} & \ldots & \frac{\partial Q_1}{\partial x_n}\\
\rule{0mm}{5mm}
\frac{\partial Q_2}{\partial x_1} & \ldots & \frac{\partial Q_2}{\partial x_n}
\end{array}\right)\]
has rank 2 for every non-zero vector $\x\in\overline{K}^n$. It follows
automatically from this that if $n\ge 4$ then the projective 
variety $Q_1=Q_2=0$ is absolutely irreducible (see Lemma \ref{arc})
and is not a cone.

The papers by Colliot-Th\'{e}l\`{e}ne, Sansuc and
Swinnerton-Dyer suggest a line of attack for the 8 variable case (see
\cite[Remark 10.5.3]{CTSSD2}). The main obstacle to this plan, 
which they failed to handle, was a purely local problem concerning
the forms $Q_1$ and $Q_2$. This we now resolve in the following
theorem.

\begin{theorem}\label{T2}
Let $k_v$ be the completion of a number field $k$ at a finite place $v$.
Let $F_v$ be the residue field of $k_v$ and assume that $\# F_v\ge 32$.
Suppose that $Q_1(X_1,\dots,X_8)$ and 
$Q_2(X_1,\dots,X_8)$ are quadratic forms over $k_v$ such that the
projective variety $Q_1(\b{X})=Q_2(\b{X})=0$ is nonsingular, and
assume further that the forms have a non-trivial common zero over $k_v$.
Then there is a form $Q$ in the pencil $a_1Q_1+a_2Q_2$ (with
$a_1,a_2\in k_v$) containing (at least) 3
hyperbolic planes. 
\end{theorem}

We will deduce Theorem \ref{T1} from Theorem \ref{T2}, but our
argument differs somewhat from that sketched by 
Colliot-Th\'{e}l\`{e}ne, Sansuc and Swinnerton-Dyer.  

While Theorem 2 can be appropriately extended to singular
intersections, it is unclear whether one can prove a corresponding
global statement. Thus the methods of the present paper seem
insufficient to handle a version of Theorem \ref{T1} for singular
intersection of two quadrics in $\proj^7$.

The next section will describe some basic facts and terminology from
the theory of quadratic forms, but we should stress at the outset that
much of our analysis involves the reduction of integral forms from $k_v$
down to $F_v$.  It follows that if $\chi(F_v)=2$ we are forced to 
consider quadratic
forms in characteristic 2, which may be less familiar to some readers
than the case of odd (or infinite) characteristic.

Throughout the paper we will write $k_v$ for the completion of a
number field $k$.  When $v$ is a non-archimedean valuation we will
write $\ov$ for the valuation ring of $k_v$, and $F_v$ for the
residue field.  We will tend to use upper case $Q$ and $L$
(and other letters) for quadratic and linear forms over
$k$ or $k_v$, and similarly lower case $q$ and $\ell$ (and other
letters) for forms over the finite field $F_v$.

\section{Quadratic Forms}\label{1QF}

In this section we recall some basic facts about the theory of
quadratic forms.  In the case of characteristic 2 the reader may wish
to treat some of our statements as exercises.

For any field $K$ we will think of a quadratic form
$q(X_1,\ldots,X_n)$ over $K$ as a polynomial of the shape
\begin{equation}\label{q1}
q(X_1,\ldots,X_n)=\sum_{1\le i\le j\le n}q_{ij}X_iX_j.
\end{equation}

It is usual to represent quadratic forms using coefficients $q_{ij}$
for every pair $i,j\le n$, subject to the symmetry constraint
$q_{ij}=q_{ji}$.  However this is inappropriate in characteristic 2
and we therefore depart from this convention.

We define $\rank(q)$ to be the least integer $m$ such that there
is a linear transformation $M\in \GL_n(K)$ for which $q(M\b{X})$ is a
function of $X_1,\ldots,X_m$ alone.  If $K$ is contained in a larger
field $K'$ then the rank over $K'$ will be the same as the
rank over $K$.  Moreover $\rank(q)<n$ if and only if $q(\b{X})=0$ has a
singular zero over $K$, that is to say a zero $\b{x}\not=\b{0}$ 
for which $\nabla
q(\b{x})=\b{0}$. This is exactly the case in which the projective
variety $\cl{Q}:q=0$ is singular, or equivalently, is a cone.  Indeed the
set of vertices for $\cl{Q}$ (as a cone in projective space) is a linear
space $\cl{V}$ of codimension $\rank(q)$. If $P$ is a nonsingular point on
$\cl{Q}$ with tangent hyperplane $\cl{H}$, then
$\cl{V}\subseteq\cl{H}$, and one sees that all points of $<P,\cl{V}>$
are vertices for $\cl{Q}\cap\cl{H}$ as a cone. The linear space
$<P,\cl{V}>$ has codimension $\rank(q)-2$ in $\cl{H}$.
Thus if we project
$\cl{Q}\cap\cl{H}$ to $\proj^{n-2}$, so as to produce a quadric
hypersurface, the corresponding quadratic form will have rank at most
$\rank(q)-2$ (and in fact equal to $\rank(q)-2$).  This produces the
following result
\begin{lemma}\label{add}
Let $q(X_1,\ldots,X_n)$ be a quadratic form of rank $r$, and suppose
$X_m=0$ is tangent to $q=0$ at a nonsingular point.  Then
$q(X_1,\ldots,X_{n-1},0)$ has rank $r-2$.
\end{lemma}

For a quadratic form in the shape (\ref{q1}) we may define an
$n\times n$ matrix $M(q)$ with
\beql{MQ}
M(q)_{ij}=\left\{\begin{array}{cc} q_{ij}, & i<j,\\ 2q_{ij}, & i=j,\\
q_{ji}, & j<i. \end{array}\right. 
\eeq
When $\chi(K)=2$ this will have diagonal entries equal to zero. We
note that $\nabla q(\b{X})=M(q)\b{X}$.  Moreover
if $T\in\GL_n(K)$ and $q_T(\b{X})=q(T\b{X})$ then
$M(q_T)=T^tM(q)T$, where $T^t$ denotes the transpose of $T$.

We also define the determinant $\det(q)$ as 
\[\det(q):=\det(M(q)),\]
so that $\det(q_T)=\det(T)^2\det(q)$.  Our definition, which is
designed to help in the case of characteristic 2, differs by powers of
2 from that which the reader may have expected.  It should be observed 
that $\det(q)$ will vanish identically if $\chi(K)=2$
and $n$ is odd.  In this case one may use the
``half-determinant'' of $q$ which we denote $\hd(q)$.  Although one can
define this in complete generality (see Leep and Schueller \cite{LS}) 
we shall need it only for
the case in which $K$ is the residue field $F_v$ of a completion $k_v$ for
some number field $k$, with respect to a non-archimedean valuation $v$.  
Then, if $q(\b{X})\in F_v[\b{X}]$ is
the reduction of some form $Q(\b{X})\in \ov[\b{X}]$, 
one can verify that $\tfrac{1}{2}\det(Q)\in \ov$.  We then define
\[\hd(q)=\overline{\tfrac{1}{2}\det(Q)},\]
(where $\overline{\theta}$ denotes the reduction to $F_v$ of 
$\theta\in \ov$). 
One may verify that this is indeed independent of the choice of
the lift $Q$ of $q$, and that $\hd(q_T)=\det(T)^2\hd(q)$. A quadratic
form is nonsingular if and only if its determinant (or in the case in
which $n$ is odd and $\chi(K)=2$, its half-determinant) is non-zero.

If the form $q$ has rank at most $r$ the corresponding matrix $M(q)$ will also
have rank at most $r$.  When $r<n$ it follows that all $(r+1)\times(r+1)$
submatrices of $M(q)$ are singular.  In the case in which $\chi(K)=2$
and $r$ is even we can say slightly more. Suppose 
that $q_0$ is a 
quadratic form in $r+1$ variables obtained by setting to zero $n-r-1$
of the variables in $q$. Then $q_0$ will be singular and
$\hd(q_0)=0$. Thus not only do all the $(r+1)\times(r+1)$ minors
vanish, but also the ``central'' $(r+1)\times(r+1)$ half-determinants.

We shall express these conditions in general by saying somewhat loosely that 
all the $(r+1)\times(r+1)$ minors of $q$ vanish.  Here one should bear
in mind that if $r$ is even and $\chi(K)=2$ then this 
refers to the half-determinant in the case of ``central'' minors.

Conversely, the vanishing of all the $(r+1)\times(r+1)$ minors of $q$
implies that $\rank(q)\le r$. To prove this it is enough to verify that if the
$(r+1)\times(r+1)$ minors of $q$ all vanish then so do those of $q_T$
for any elementary matrix $T$. These facts are enough to show that 
\beql{r+2}
\rank\big(q(X_1,\ldots,X_m)+X_{m+1}X_{m+2}\big)=
\rank\big(q(X_1,\ldots,X_m)\big)+2.
\eeq

If $q(\b{x})=0$ has a nonsingular zero over $K$ there is a linear
transformation $T\in \GL_n(K)$ for which $q(T\b{X})$ takes the form
$X_1X_2+q'(X_3,\ldots,X_n)$.  We then say that $q$ 
``splits off a hyperbolic plane''. If $q$ splits off $m$ hyperbolic
planes then the form $q'$ will split off $m-1$ hyperbolic planes.

When $K$ is finite, as is often the case of interest to us, the
Chevalley--Warning theorem implies that $q$ has a nonsingular zero
over $K$ as long as $\rank(q)\ge 3$.  It follows that there is a matrix
$T$ such that
\[q(T\b{X})=X_1X_2+\ldots+X_{2s-1}X_{2s}+cX_{2s+1}^2\]
for some $c\not=0$, if $q$ has odd rank equal to $2s+1$, or such
that $q(T\b{X})$ takes one of the shapes
\[q(T\b{X})=X_1X_2+\ldots+X_{2s-1}X_{2s}\]
or
\[q(T\b{X})=X_1X_2+\ldots+X_{2s-3}X_{2s-2}+n(X_{2s-1},X_{2s})\]
if $\rank(q)=2s$ is even.  Here $n(X,Y)$ is an anisotropic form.  We
note that, when $K$ is finite, if $n(X,Y)$ and $n'(X,Y)$ are two
anisotropic forms there will be a linear transform in $GL_2(K)$ taking
$n$ to $n'$.  Moreover there is a linear transform in $GL_4(K)$ taking
$n(X_1,X_2)+n'(X_3,X_4)$ to $X_1X_2+X_3X_4$.  We shall be somewhat lax
in our notation for anisotropic forms, writing $n(X,Y)$ for a generic
form of this type, not necessarily the same at each occurrence.
Moreover we shall sometimes use the notation $N(X,Y)$ for a binary form over
$\ov$ whose reduction to $F_v$ is anisotropic.

The following result will be useful in recognizing when a form over
$\ov$ splits off hyperbolic planes.
\begin{lemma}\label{lift}
Let $\pi$ be a uniformizing element for $k_v$ and suppose that the
quadratic form
$Q(X_1\ldots,X_n)\in \ov[\b{X}]$ satisfies
\begin{eqnarray*}
\lefteqn{Q(\b{X})\equiv
X_1X_2+\ldots+X_{2s-1}X_{2s}+\widetilde{Q}(X_{2s+1},\ldots,X_n)}\hspace{2cm}\\
&&+\pi\sum_{i=1}^{2s}X_iL_i(X_1,\ldots,X_n)\modd{\pi^2}
\end{eqnarray*}
for some quadratic form $\widetilde{Q}$ over $\ov$.  Then there exists
$T\in\GL_n(\ov)$ such that
\[Q(T\b{X})=X_1X_2+\ldots+X_{2s-1}X_{2s}+Q_0(X_{2s+1},\ldots,X_n)\]
with $Q_0\equiv \widetilde{Q}\modd{\pi^2}$.

In particular, if $\overline{Q}$ splits off at least $s$
hyperbolic planes over $F_v$, then so does $Q$ over $k_v$.  Indeed if
$\overline{Q}$ has rank at least 7 then $Q$ splits off at least 3
hyperbolic planes. 

\end{lemma}
Unfortunately if $Q_1$ and $Q_2$ are forms in 8 variables
with coefficients
in $\ov$ it is possible that every linear combination
$a\overline{Q_1}+b\overline{Q_2}$ has rank at most 6, even if the
variety $Q_1=Q_2=0$ is nonsingular over $k_v$.  This is where the 
difficulty in proving Theorem \ref{T2} lies.

We may establish Lemma \ref{lift} by adapting the argument for 
Hensel's Lemma.  We show
inductively that for every positive integer $h$ there is a 
$T_h\in\GL_n(\ov)$ and linear forms $L_i^{(h)}(X_1,\ldots,X_n)$
over $\ov$ such that
\begin{eqnarray*}
\lefteqn{Q(T_h\b{X})\equiv
X_1X_2+\ldots+X_{2s-1}X_{2s}+Q_h(X_{2s+1},\ldots,X_n)}\hspace{2cm}\\
&&+\pi^h\sum_{i=1}^{2s}X_iL_i^{(h)}(X_1,\ldots,X_n)\modd{\pi^{h+1}}
\end{eqnarray*}
with $Q_h\equiv \widetilde{Q}\modd{\pi^2}$. Once this is established one may
choose a convergent subsequence from the $T_h$ and the lemma will
follow.

To prove the claim we observe that the case $h=1$ is
immediate. Generally if $U_h\in\GL_n(\ov)$ corresponds to the
substitution
\begin{eqnarray*}
X_{2i-1}&\rightarrow &
X_{2i-1}+\pi^hL_{2i}^{(h)}(X_1,\ldots,X_n),
\;\;\; (1\le i\le s)\\
X_{2i}&\rightarrow &
X_{2i}+\pi^hL_{2i-1}^{(h)}(X_1,\ldots,X_n),
\;\;\; (1\le i\le s)\\
X_i&\rightarrow& X_i,\;\;\; (2s<i\le n),
\end{eqnarray*}
then it suffices to take $T_{h+1}=U_hT_h$.  
\bigskip

There is one further result which will be useful in finding forms which
split off three hyperbolic planes.
\begin{lemma}\label{nsq}
Let $Q(X_1,\ldots,X_8)$ be a quadratic form over $k_v$ whose
determinant is not a square in $k_v$.  Then $Q$ 
contains at least three hyperbolic planes.
\end{lemma}

The form $Q$ must be nonsingular. However any nonsingular form in 8 
variables over $k_v$ will split off two
hyperbolic planes, leaving a form in 4 variables, so that
$Q=X_1X_2+X_3X_4+Q'(X_5,\ldots,X_8)$.  Then $\det(Q')=\det(Q)$ is a
non-square, which implies that $Q'$ is isotropic, and hence splits off
a further hyperbolic plane.
\bigskip

In view of Lemma \ref{nsq}, it would suffice for the proof of Theorem
\ref{T2} to find any form $aQ_1+bQ_2$ whose determinant is not a
square.  Unfortunately it is possible that all forms in the pencil
have square determinant, as the example
\[Q_1=X_1^2-X_2^2+X_3^2-4X_4+X_5^2-7X_6^2+X_7^2-10X_8^2\]
\[Q_2=X_1X_2+X_3X_4+X_5X_6+X_7X_8\]
over $\mathbb{Q}_3$ shows. These forms have a common zero at
$(2,1,0,0,2,-1,0,0)$, as in the hypotheses of Theorem \ref{T2}.

As we have remarked, our proof of Theorem \ref{T2} will involve the
reductions $\overline{Q}$ of quadratic forms $Q$, defined over
$\ov$. The argument will require a great many invertible linear
transformations of variables. In this context we will use the fact
that any element of $\GL_n(F_v)$ can be lifted to $\GL_n(\ov)$.  Thus
whenever we refer to a ``change of variables'' or a ``substitution''
among the variables of such a form over $F_v$, we mean that one applies an
appropriate element of $\GL_n(\ov)$ which reduces to the relevant
mapping over $F_v$.

Since we will use a succession of changes of variables it will become
cumbersome to use different notation for all the different forms that
arise.  We will therefore abuse notation by saying for example that
the change of variables $X_1\rightarrow X_1+X_2$ transforms
$q(X_1,X_2)=X_1^2$ into $q(X_1,X_2)=X_1^2+2X_1X_2X_2^2$, rather than
producing $q'(X_1,X_2)=X_1^2+2X_1X_2X_2^2$, say. Thus we may have a
number of forms, all denoted by ``$q$'', which are not actually the
same.  We trust that this will not cause confusion.

\section{Pairs of Quadratic Forms}\label{PQF}

Given two quadratic forms $q_1(X_1,\ldots,X_n),q_2(X_1,\ldots,X_n)$
defined over a field $K$ we may consider the pencil
$\cP=<q_1,q_2>=<q_1,q_2>_K$, which is the set of all linear
combinations $aq_1+bq_2$ with $a,b\in K$, not both zero.
We will also consider $\cP^*=<q_1,q_2>_{\overline{K}}$, where 
$\overline{K}$ is the algebraic closure of $K$.
There are three notions of rank that we attach to such a pencil. We
define $R(\cP)$ as the least integer $m$ such that there
is a linear transformation $T\in \GL_n(K)$ for which $q(T\b{X})$ is a
function of $X_1,\ldots,X_m$ alone, for every $q\in\cP^*$. Thus $m$ is
the least integer such that there are variables $X_1,\ldots,X_m$ which
suffice to represent every form in the pencil.  We define a second
number $r(\cP)$ as $\max \rank(q)$, where $q$ runs over all forms in
$\cP^*$.  The third notion of rank which we shall use is
$r_{\min}(\cP)$ defined similarly as $\min\rank(q)$, where $q$ runs 
over all forms in $\cP^*$. 

For a pencil $\cP$ generated by $q_1$ and $q_2$ we will
write $R(q_1,q_1)=R(\cP)$, $r(q_1,q_1)=r(\cP)$ and $r_{\min}(q_1,q_1)
=r_{\min}(\cP)$.
Clearly one has 
\[r_{\min}(\cP)\le r(\cP)\le R(\cP)\]
in every case.  We remark
that if $\rank(tq_1+q_2)\le r$ for every value of $t$ then the
$(r+1)\times(r+1)$ minors of $tq_1+q_2$ all vanish identically as
polynomials in $t$, whence the minors of $aq_1+bq_2$ also vanish
identically, yielding $r(q_1,q_2)\le r$.

As an example of these notions of rank one may consider the forms
$q_1(X_1,X_2,X_3)=X_1X_2$ and $q_2(X_1,X_2,X_3)=X_1X_3$, which
generate a pencil $\cP$ with $r_{\min}(\cP)=r(\cP)=2$ and $R(\cP)=3$. 

When the variety $q_1=q_2=0$ is nonsingular we can describe these
ranks precisely.
\begin{lemma}\label{LBP}
Suppose that $Q_1(X_1,\ldots,X_n)$ and $Q_2(X_1,\ldots,X_n)$ are quadratic
forms over a field $K$ of characteristic not equal to 2. Define the
binary form
\beql{Fdef}
F(x,y)=F(x,y;Q_1,Q_2)=\det\big(xM(Q_1)+yM(Q_2)\big)
\eeq
of degree $n$.  Then if the variety $Q_1=Q_2=0$ is nonsingular the
form $F$ does not vanish identically, and has distinct linear factors
over the algebraic completion $\overline{K}$.  Moreover $r(Q_1,Q_2)=n$
and $r_{\min}(Q_1,Q_2)=n-1$.
\end{lemma}
This follows from Heath-Brown and Pierce \cite[Proposition 2.1]{HBP}, 
for example.
\bigskip

We now give a general condition for an intersection of quadrics to be
absolutely irreducible of codimension 2.
\begin{lemma}\label{arc}
Let $Q_1(X_1,\ldots,X_n)$ and $Q_2(X_1,\ldots,X_n)$ be two quadratic forms
over an algebraically closed field $K$, such that 
\[r_{\min}(Q_1,Q_2)\ge 3\;\;\;\mbox{and}\;\;\; r(Q_1,Q_2)\ge 5.\]
Then the projective variety 
\[\cl{V}: Q_1=Q_2=0\]
is absolutely irreducible of codimension
2. Moreover if $n\ge 4$ and $\cl{V}$ is nonsingular then it is 
absolutely irreducible of codimension 2.
\end{lemma}

If $\cl{V}$ fails to be absolutely irreducible of codimension 2 then it must
contain either a quadric of codimension 2, or a linear space of 
codimension 2.   Suppose firstly that $\cl{V}$ contains a quadric of 
codimension 2, given by the simultaneous vanishing of $Q(\b{X})$ and 
$L(\b{X})$ say.
Then we may write $Q_i=c_iQ+LM_i$ for suitable constants $c_i$ and
linear forms $M_i$. We now set $a_i=c_i$ for $i=1,2$, unless
$c_1=c_2=0$, in which case we take $a_1=1$ and $a_2=0$.  Then $a_1Q_1-a_2Q_2$
is a multiple of $L$ and hence has rank at most 2, 
contradicting our hypotheses.  Now consider the second case, 
in which $\cl{V}$ contains a linear space of codimension 2,
given by the simultaneous vanishing of
$L_1(\b{X})$ and $L_2(\b{X})$ say. Then $Q_1$ and
$Q_2$ must take the shape $Q_i=L_1M_{i1}+L_2M_{i2}$ for $i=1$ and 2,
where $M_{ij}$ are suitable linear forms. Hence
every quadratic in the pencil $<Q_1,Q_2>$ will have rank at 
most 4.  This contradicts our hypothesis that $r(Q_1,Q_2)\ge 5$.

Finally, if $\cl{V}$ is nonsingular and $n\ge 4$ we know
from Lemma \ref{LBP} that $r_{\min}(Q_1,Q_2)=n-1\ge 3$,
so that $\cl{V}$ cannot contain a quadric of codimension 2. Moreover 
if $\cl{V}$
contained the linear space $L_1=L_2=0$, so that we can write
$Q_i=L_1M_{i1}+L_2M_{i2}$ for $i=1$ and 2, there would be singular
points wherever 
\[L_1=L_2=M_{11}M_{22}-M_{12}M_{21}=0.\]
This contradiction completes the proof of Lemma~\ref{arc}.
\bigskip

It will be
important for us to understand the structure of pencils for which
$r(\cP)<R(\cP)$.  In this context we have the following result.

\begin{lemma}\label{lemshape1}
Suppose that
$\cP$ is a pencil of quadratic forms in $n$ variables
over a field $F$, and assume that $\# F\ge n$. Write $r(\cP)=r$ and
$R(\cP)=R$ and suppose that $r<R$.  Assume further either 
that $\chi(F)\not=2$, or that $r$
is even and $F$ is the residue field of the
completion of a number field under a non-archimedean
valuation.  Then we can choose a basis for $\cP$ such that
$\rank(q_1)=r$.  For any such basis $q_1,q_2$ there is
an invertible change of variables over $F$ so that 
\[q_1(X_1,\ldots,X_n)=q_1'(X_1,\ldots,X_r)\]
and
\[q_2(X_1,\ldots,X_n)=q_2'(X_1,\ldots,X_{R-1})+X_rX_R\] 
for certain quadratic forms $q_1',q_2'$ defined over $F$.  
\end{lemma}

If $q_1$
and $q_2$ generate $\cP$ and $r(q_1,q_2)=r$, the $r\times r$ minors of a linear 
combination $q_1+xq_2$ cannot all vanish identically in $x$. 
Since at least one of these minors is a non-zero polynomial of degree at
most $r$ in $x$
there can be at most $r$ values of $x$ for which
$\rank(q_1+xq_2)<r$.  Thus if $\# F\ge n\ge R>r$ there will be
a linear combination $q_1'=q_1+xq_2$ defined over $F$ with rank
$r$.  Thus it will be possible to choose a basis $q_1',q_2$ for $\cP$
in which $\rank(q_1')=r$.

Now, assuming that $\rank(q_1)=r$, we make an appropriate 
change of variables so that only $X_1,\ldots,X_R$ appear in the forms
in $\cP^*$. After a further
change of variables we can then write
$q_1=q_1'(X_1,\dots,X_r)$ and 
\beql{q2f}
q_2=q_3(X_1,\dots,X_r)+\sum_{i=r+1}^RX_i\ell_i(X_1,\dots,X_r)
+q_4(X_{r+1},\ldots,X_R)
\eeq
say, with appropriate linear forms $\ell_i$. With the 
notation (\ref{MQ}) the matrix for $q_1+xq_2$ then takes
the form
\beql{mat}
\left(\begin{array}{c|c} \rule[-2mm]{0mm}{1mm}M(q_1'+xq_3) & x\ell\\ \hline
\rule{0mm}{5mm}x\ell^t & xM(q_4)\end{array}\right).
\eeq
Since $q_1+xq_2$ has rank at most $r$ for every value of $x$ we deduce
that all the $(r+1)\times(r+1)$ minors vanish identically in $x$
(with the obvious interpretation for the central minors when $r$ is even
and $\chi(F)=2$).  Suppose now that $q_4$ contains a term $cX_iX_j$
with $r<i\le j\le R$, so that $c=(q_4)_{ij}$.  If $i\not=j$, or if
$\chi(D)\not=2$,  we consider the minor
formed from rows $1$ to $r$ and $i$, and columns 1 to $r$ and $j$.
This will be a polynomial in $x$, in which the coefficient of $x$ is
$c\det(q_1')$.  Since the polynomial vanishes identically we deduce
that $c\det(q_1')=0$.  Similarly if $i=j$ and $\chi(F)=2$ (so that $r$ is
even) the half-determinant corresponding to the $r+1$ variables
$X_1,\ldots,X_r$ and $X_i$ must vanish.  However this half determinant
is a polynomial in $x$ with linear term $c\det(q_1')x$, and again we 
deduce that $c\det(q_1')=0$.
However since $q_1$ has rank $r$, and we have avoided the
case in which $\chi(F)=2$ and $r$ is odd, we will have
$\det(q_1')\not=0$, so that $c=0$.  We therefore deduce that all
coefficients of $q_4$ vanish.

We see now that (\ref{q2f}) simplifies to
\[q_2=q_3(X_1,\dots,X_r)+\sum_{i=r+1}^RX_i\ell_i(X_1,\dots,X_r).\]
If all the forms $\ell_{r+1},\ldots, \ell_R$ were identically zero
we would have $R(\cP)\le r$, contrary to our assumption.  We therefore
suppose that $\ell_R$, say, is not identically zero.  Without loss of
generality we may assume that $\ell_R$ involves $X_r$ with non-zero
coefficient.  Indeed, after a change of variable among
$X_1,\ldots,X_r$ we can then assume that $\ell_R(X_1,\ldots,X_r)=X_r$.
It then follows that $q_1$ and $q_2$ are of the shape given in Lemma
\ref{lemshape1}.
\bigskip

For the rest of the present section, the pairs of forms
$q_1,q_2$ will be defined over a field 
$F_v$, which will be the residue field
of the completion $k_v$ of a number field $k_v$ with respect to a
non-archimedean valuation.  Some of the results we prove will in fact
be valid in a more general setting, but this will suffice for our needs.
\bigskip

We now examine the remaining case, in which $\chi(F_v)=2$ and $r$ is
odd. By extending the argument above we will prove the following
structure result.
\begin{lemma}\label{lemshape2}
Let $F_v$ and $\cP$ be as in Lemma \ref{lemshape1}, except that
$\chi(F_v)=2$ and that $r$ is odd.  Then there is a basis $<q_1,q_2>$
for $\cP$ in which $\rank(q_1)=r$.  For any such basis there is
an invertible change of variables over $F_v$ such that 
\[q_1(X_1,\ldots,X_n)=q_1'(X_1,\ldots,X_r)\]
and $q_2$ is in one of the shapes
\beql{q2sh}
q_2(X_1,\ldots,X_n)=q_2'(X_1,\ldots,X_{R-1})+X_rX_R,
\eeq
or
\beql{q2sh1}
q_2(X_1,\ldots,X_n)=q_2'(X_1,\ldots,X_r)+X_{r+1}^2,
\eeq
or
\beql{q2sh2}
q_2(X_1,\ldots,X_n)=q_2'(X_1,\ldots,X_r)+X_{r+1}^2+X_rX_{r+1}.
\eeq
When $r\le R-2$ we can always take $q_2$ to be of the shape (\ref{q2sh}).
\end{lemma}

We start the proof by choosing generators $q_1$ and $q_2$ for $\cP$ as
before, and expressing $q_2$ in the shape (\ref{q2f}).  As before the
$(r+1)\times(r+1)$ minors of (\ref{mat}) must all vanish identically
in $x$.  We begin by
considering an ``off-diagonal'' term, $cX_iX_j$ say, in
$q_4(X_{r+1},\ldots,X_R)$.  Thus we will assume that $r+1\le i<j\le
R$.  Choose any
sets $I_0,J_0\subset\{1,\ldots,r\}$ of cardinality $r-1$ and let 
$I=I_0\cap\{i,j\}$ and $J=J_0\cap\{i,j\}$. If $I_0\not=J_0$ then the
$I,J$ minor of (\ref{mat}) is a determinant which is a polynomial in
$x$, which must vanish identically.  The coefficient of $x^2$ will 
be $-c^2$ times the $I_0,J_0$ minor of $q_1'$, so that this product
must vanish.  Similarly if $I_0=J_0$ the $I,I$ minor of (\ref{mat})
is a half-determinant which is again a polynomial in
$x$.  The coefficient of $x^2$ will similarly be $-c^2$ times the $I_0,I_0$
minor of $q_1'$, which must also vanish. If all the $(r-1)\times(r-1)$ minors
of $q_1'$ were to vanish we would have $\rank(q_1')\le r-2$, contrary to
hypothesis. We therefore conclude that
$c=0$.  Thus $q_4$ will be a diagonal form.  Since $F_v$ is a finite
field of characteristic 2 there is therefore a change of
variable reducing $q_4$ to the shape $cX_{r+1}^2$. Now, if $c=0$, or
if any of the forms
$\ell_{r+2},\ldots,\ell_R$ is not identically zero, we can proceed as in the
proof of Lemma \ref{lemshape1}, so as to put $q_2$ into the form
(\ref{q2sh}).  Moreover if $\ell_{r+2},\ldots,\ell_R$ all vanish identically
we see that we must have $R=r+1$.

Thus it remains to consider the case in which $R=r+1$ and
\[q_2(X_1,\ldots,X_n)=q_2'(X_1,\ldots,X_r)+cX_R^2+X_R\ell(X_1,\ldots,X_r)\]
with $c\not=0$. Since every element of $F_v$ is a square we may
replace $X_R$ by $c^{-1/2}X_R$ so as to reduce to
the case $c=1$.  Then, if $\ell$ vanishes identically we obtain a form of type
(\ref{q2sh1}), while if $\ell$ does not vanish identically we can make a
change of variables so as to replace $\ell(X_1,\dots,X_r)$ by $X_r$,
giving us a form of type (\ref{q2sh2}). This completes the proof of
Lemma \ref{lemshape2}.
\bigskip

We can say a little more about the structure of our forms when $q_2$
has the shape (\ref{q2sh}).
\begin{lemma}\label{r-2}
Let $\cP$ be a pencil with $r(\cP)=r<R=R(\cP)$, and suppose either that
$r$ is even or that $r\le R-2$.  Then $r\ge 2$.  Moreover if $\rank(q_1)=r$
we can make a change of variable so that
\[q_1=q_3(X_1,\ldots,X_{R-2})+X_{R-1}\ell(X_1,\ldots,X_{R-1})\]
and
\[q_2=q_4(X_1,\ldots,X_{R-2})+X_{R-1}X_R\]
with $\rank(q_3)=r(q_3,q_4)=r-2$ and $R(q_3,q_4)=R-3$ or $R-2$.  In
particular if $r=2$ the forms $q_3$ and 
$q_4$ vanish identically.
\end{lemma}

We begin the proof by observing that, according to Lemmas
\ref{lemshape1} and \ref{lemshape2}, we can put $q_2$ into the shape  
(\ref{q2sh}).  We may then rewrite the forms as
\[q_1=q_1''(X_1,\ldots,X_{r-1})+X_r\ell_1(X_1,\ldots,X_r)\]
and
\[q_2=q_2''(X_1,\ldots,X_{r-1},X_{r+1},\ldots,X_{R-1})
+X_r\big(X_R+\ell_2(X_1,\ldots,X_{R-1})\big).\]
We then replace $X_R$ by $X_R+\ell_2(X_1,\ldots,X_{R-1})$ and re-number the
variables so as to interchange $X_r$ and $X_{R-1}$.  This puts $q_1$ 
and $q_2$ into the shape given in the lemma. For any $x$ we then have
\[xq_1+q_2=xq_3+q_4+X_{R-1}\big(x\ell(X_1,\ldots,X_{R-1})+X_R\big).\]
Replacing $X_R$ by $X_R+x\ell(X_1,\ldots,X_{R-1})$ we see from (\ref{r+2}) that
\[\rank(xq_1+q_2)=\rank(xq_3+q_4)+2\]
for all $x$, whence we must have $r(q_1,q_2)\ge 2$ and
$r(q_3,q_4)=r-2$, as required.  Moreover it is evident that
$r=\rank(q_1)\le \rank(q_3)+2$, so that $\rank(q_3)=r-2$.  Clearly we have
$R(q_3,q_4)\le R-2$, and if it were possible to write $q_3$ and $q_4$
using only $R-4$ variables then the expressions for $q_1$ and $q_2$
given in the lemma would allow us to write $q_1$ and $q_2$ with only
$R-1$ variables.  Since this is impossible we conclude that
$R(q_3,q_4)\ge R-3$.
\bigskip

One final lemma belongs in this section.
\begin{lemma}\label{sz}
Let $q_1$ and $q_2$ be  quadratic forms in $4$ variables, defined over a field
$F_v$ as above.  Then if
$r(q_1,q_2)<4$ the forms will have a singular common zero
over $F_v$; that is to say there is a non-zero vector $\x\in F_v^4$
with $q_1(\x)=q_2(\x)=0$ and
$\nabla q_1(\x)$ and $\nabla q_2(\x)$ proportional. 
\end{lemma}
It should be observed that when $\chi(F_v)=2$ this result may fail for
forms in $n\not=4$ variables.  For example, when $n=2$ the forms
$q_1=X_1^2$ and $q_2=X_2^2$ have $r=1<n$, but have no singular common
zero over $F_v$ (indeed they have no non-trivial common zero).
Equally, when $n=6$ the forms
\[q_1(\b{X})=X_1X_2+X_3X_4+X_5^2\]
and
\[q_2(\b{X})=X_1X_2+cX_1X_3+X_2X_4+X_6^2\]
have no singular common zero over $F_v$ provided that $c\in F_v$ is
chosen so that the polynomial $T^2+T+c$ is irreducible.

We can choose
$\x=(0,0,0,1)$ when $q_2$ is of the shape (\ref{q2sh}).  
According to Lemmas \ref{lemshape1} and \ref{lemshape2} 
it therefore remains to consider the case in which
$\chi(F_v)=2$ and $r=3$, and $q_2$ takes one of the forms (\ref{q2sh1})
or (\ref{q2sh2}).   

We begin by examining the first case, in which
\[q_1(X_1,\ldots,X_4)=q_1'(X_1,X_2,X_3)\;\;\mbox{and}\;\;
q_2(X_1,\ldots,X_4)=q_2'(X_1,X_2,X_3)+X_4^2.  \]
We begin by changing variables so that 
\[q_2'(X_1,X_2,X_3)=q_3(X_1,X_2)+cX_3^2.\]
We can now write $c=d^2$ and replace $X_4$ by $X_4+dX_3$ to get
\[q_2(X_1,\ldots,X_4)=q_3(X_1,X_2)+X_4^2.\]
Then $(0,0,1,0)$ will be a
singular common zero if $q_1'(0,0,1)=0$.  Otherwise we may take
\[q_1(X_1,\ldots,X_4)=q_4(X_1,X_2)+\ell(X_1,X_2)X_3+eX_3^2,\]
say, with $e\not=0$.  We now choose $x_1,x_2\in F_v$, not both zero, such that
$\ell(x_1,x_2)=0$. There then exist $x_3,x_4\in F_v$ with
$ex_3^2=q_4(x_1,x_2)$ and $x_4^2=q_3(x_1,x_2)$, so that
$q_1(\x)=q_2(\x)=0$.  However we have arranged that the third and
fourth entries of $\nabla q_i(\x)$ vanish for $i=1$ and $2$.  Moreover
$\nabla q(\b{X})$ is automatically orthogonal to $\b{X}$, for any
quadratic form $q$ over a field of characteristic 2.  Hence
\[x_1\frac{\partial q_1(\x)}{\partial x_1}+
x_2\frac{\partial q_1(\x)}{\partial x_2}=0,\;\;\;\mbox{and}\;\;\;
x_1\frac{\partial q_2(\x)}{\partial x_1}+
x_2\frac{\partial q_2(\x)}{\partial x_2}=0.\]
This is enough to ensure that $\nabla q_1(\x)$ and $\nabla q_2(\x)$
are proportional.

In the case of (\ref{q2sh2}) we write
\[q_1(\b{X})=q_3(X_1,X_2)+X_3\ell_1(X_1,X_2,X_3)\]
and
\[q_2(\b{X})=q_4(X_1,X_2)+X_3\big(\ell_2(X_1,X_2,X_3)+X_4\big)+X_4^2,\]
whence
$\det(q_2+xq_1)=x^2\det(q_4+tq_3)$.  However $\det(q_2+xq_1)$ vanishes
identically since $r(q_1,q_2)<4$, and we therefore see that
$\det(q_4+xq_3)$ also vanishes identically, yielding $r(q_3,q_4)\le 1$. Since
$\chi(F_v)=2$, we know that every element of $F_v$ is a square, and it
follows that we can write
\[q_1(\b{X})=\ell_3(X_1,X_2)^2+X_3\ell_1(X_1,X_2,X_3)\]
and
\begin{eqnarray*}
q_2(\b{X})&=&\ell_4(X_1,X_2)^2+X_3(\ell_2(X_1,X_2,X_3)+X_4)+X_4^2\\
&=&(\ell_4(X_1,X_2)+X_4)^2+X_3(\ell_2(X_1,X_2,X_3)+X_4).
\end{eqnarray*}
We may then complete the proof of the lemma by choosing a non-zero 
vector $\x\in F_v^4$ such that $x_3=\ell_3=\ell_4+x_4=0$ .

\section{$v$-Adically Minimized Pairs of Forms}\label{secmin}

Since our analysis is based on reduction to $F_v$ we begin by
multiplying the forms $Q_1,Q_2$ in Theorem \ref{T2}
by a suitable scalar so as to give them
coefficients in $\ov$. For such integral forms we will write
$q_1=\overline{Q_1}$ and $q_2=\overline{Q_2}$ for the reductions to
$F_v$. We retain this notation through to the end of Section \ref{secR7}.

With these conventions we can now explain the fundamental idea
behind our proof.  We begin by choosing a suitable model for the pencil
$<Q_1,Q_2>$ over $\ok$.  For example, one may divide $Q_1$ by a
suitable power of $\pi$ so as to reach a form for which $q_1$ does not
vanish identically. This choice of model will depend on a $v$-adic
minimization technique from the work of Birch, Lewis and Murphy \cite{BLM}. We
then consider a large number of cases, depending principally on the
values of $r=r(q_1,q_2)$ and $R=R(q_1,q_2)$. If our model is suitably
chosen it turns out that small values of $r$ and $R$ cannot occur.  In
the remaining cases we prove that some form in the pencil $<Q_1,Q_2>$
splits off three hyperbolic planes.  Often we will be able to do this
by showing that some form $aq_1+bq_2$ over $F_v$ splits off three
hyperbolic planes, so that we can apply Lemma \ref{lift}. However other
cases will require more work.

There is one particularly easy case. If $r(q_1,q_2)\ge 7$, and
assuming that $\# F_v\ge 8$, there will be a form $aq_1+bq_2$ with
$a,b\in F_v$ whose rank is at least 7.  Then we 
can conclude from Lemma \ref{lift} and that $Q$ splits off 3
hyperbolic planes, as required for Theorem \ref{T2}. 

We record this observation
as follows.
\begin{lemma}\label{r7}
If $Q_1$ and $Q_2$ are defined over $\ov$ and have $r(q_1,q_2)\ge 7$,
then there is a form in the pencil $<Q_1,Q_2>$ which
splits off three hyperbolic planes.
\end{lemma}
\bigskip

We now describe the $v$-adic minimization 
process given by Birch, Lewis and Murphy \cite{BLM}.  We 
start with two quadratic forms $Q_1(X_1,\ldots,X_n)$ and
$Q_2(X_1,\ldots,X_n)$ defined over $k_v$, the completion of a number
field $k$ with respect to a non-archimedean valuation $v$. We will
assume that the variety $Q_1=Q_2=0$ is nonsingular, in the sense
specified in Section \ref{sec-intr}.  We will retain this hypothesis
throughout the paper, without further comment.

Since the variety $Q_1=Q_2=0$ is nonsingular, Lemma \ref{LBP} tells us
that $F(x,y)$ does not vanish identically, and has no
repeated factors.  We now set 
\beql{Ddef}
\Delta(Q_1,Q_2)={\rm Disc}(F(x,y;Q_1,Q_2))
\eeq
which will be a non-zero element of $k_v$. Indeed if the forms $Q_1$
and $Q_2$ are defined over $\ov$ then $\Delta(Q_1,Q_2)$ will also be in
$\ov$.

For any matrices $U\in \GL_2(k_v)$ and $T\in\GL_n(k_v)$ we define
actions on pairs of quadratic forms $Q_1,Q_2$ by setting
\[(Q_1,Q_2)^U=(U_{11}Q_1+U_{12}Q_2,U_{21}Q_1+U_{22}Q_2)\]
and 
\[(Q_1(\b{X}),Q_2(\b{X}))_T=(Q_1(T\b{X}),Q_2(T\b{X})).\]
Notice that $Q_1=Q_2=0$ defines a smooth variety with a non-trivial
point over $k_v$ if and only if the same is true for the forms
$(Q_1,Q_2)^U_T$.  Similarly the pencil defined over $k_v$ by $Q_1,Q_2$ contains
a form which splits off 3 hyperbolic planes if and only if the same is true
for the quadratics  $(Q_1,Q_2)^U_T$. Thus we may transform our pair
$Q_1,Q_2$ in this way in the hope of producing forms of a convenient
shape. We shall say that a pair of integral forms $Q_1,Q_2$ is ``minimized''
if there are no transforms $U$ and $T$ such that $(Q_1,Q_2)^U_T$ is
also integral and
\[|\Delta((Q_1,Q_2)^U_T)|_v>|\Delta(Q_1,Q_2)|_v.\]
It is clear that there is always a pair of transforms $U$ and $T$ such that
$(Q_1,Q_2)_T^U$ is minimized.  Indeed it may happen that there are
many quite different pairs $U,T$ which may be used.
Since one can compute in general that
\[\Delta((Q_1,Q_2)^U_T)=(\det(U))^{n(n-1)}(\det(T))^{4(n-1)}\Delta(Q_1,Q_2),\]
the condition for a pair of integral forms $Q_1,Q_2$ to be minimized
is that there are no matrices $U,T$ for which $(Q_1,Q_2)^U_T$ are
integral and such that 
\[|\det(U)|_v^n|\det(T)|_v^4>1.\]
We observe for future reference that if 
\beql{TU=}
|\det(U)|_v^n|\det(T)|_v^4=1
\eeq
and $Q_1,Q_2$ is a minimized pair of integral forms, then $(Q_1,Q_2)_T^U$
will also be minimized, provided of course that the resulting 
forms are integral.
\bigskip

From now on we will restrict to the case $n=8$, for which the above
condition says that there are no suitable $T,U$ with
\beql{mincon}
|\det(U)|_v^2|\det(T)|_v>1.
\eeq
Continuing with the
notation above we have the following simple criteria for a pair of
forms not to be minimized.

\begin{lemma}\label{notmin}
Suppose that the projective variety $q_1=q_2=0$ contains a linear space
of projective dimension at least 4 defined over $F_v$.  Then the 
pair $Q_1,Q_2$ is not minimized. In particular, 
if $R(q_1,q_2)\le 3$ then the pair $Q_1,Q_2$ is not minimized.
\end{lemma}

A convenient condition equivalent to the existence of a linear space
of projective dimension at least 4 is that the forms $q_1$ and $q_2$
can be written as
\beql{3L}
\begin{array}{ccc}
q_1(\b{X})&=&\ell_1(\b{X})\lambda_1(\b{X})+\ell_2(\b{X})\lambda_2(\b{X})
+\ell_3(\b{X})\lambda_3(\b{X})\\
\rule{0mm}{5mm}
q_2(\b{X})&=&\ell_1(\b{X})\mu_1(\b{X})+\ell_2(\b{X})\mu_2(\b{X})
+\ell_3(\b{X})\mu_3(\b{X})
\end{array}
\eeq
for suitable linear forms $\ell_i,\lambda_i,\mu_i$ defined over $F_v$. 

To prove the lemma we assume
for simplicity that $\ell_1,\ell_2$ and $\ell_3$ are linearly independent.  In
the alternative case we will have a similar argument involving fewer
linear forms.  We apply a transform $T_1\in\GL_8(\ov)$ so as to
replace $\ell_i(\b{X})$ by $X_i$ for $1\le i\le 3$.  If $\pi$ is a
uniformizing element for $k_v$ we can then write
\[Q_1(T_1\b{X})=X_1L_1(\b{X})+X_2L_2(\b{X})+X_3L_3(\b{X})+\pi Q_3(\b{X})\]
for suitable linear forms $L_i$ and a quadratic form $Q_3$, all
defined over $\ov$.  There is also an analogous expression for $Q_2$.
We now define 
\[T_2=\diag(\pi,\pi,\pi,1,1,1,1,1)\]
and $T=T_2T_1$, so
that $|\det(T)|_v=|\pi|_v^3.$  We then see that both $Q_1(T\b{X})$ and
$Q_2(T\b{X})$ are divisible by $\pi$.  Hence if we take
$U=\diag(\pi^{-1},\pi^{-1})$ then $(Q_1,Q_2)^U_T$ is a pair of integral
forms.  However since $|\det(U)|_v=|\pi|_v^{-2}$ the
condition (\ref{mincon}) is satisfied.  Thus the pair $Q_1,Q_2$ is not
minimized.

When $R(q_1,q_2)\le 3$ we can write $q_1$ and $q_2$ in terms of just 3
variables $X_1,X_2,X_3$ so that the 4-plane $X_1=X_2=X_3=0$ is
contained in the variety $q_1=q_2=0$.  Lemma \ref{notmin} then follows.
\bigskip

There is a further instance in which one can see that a pair of forms
is not minimized, given by the following lemma.
\begin{lemma}\label{x8}
Suppose that $R(q_1,q_2)=R\le 7$ and that the forms $Q_1,Q_2$ take the shape
\begin{eqnarray}\label{x8shape}
Q_i(X_1,\ldots,X_8)&=&G_i(X_1,\ldots,X_R)
+\pi\sum_{j=1}^RX_jL_j^{(i)}(X_{R+1},\ldots,X_8)\nonumber\\
&&\hspace{1cm}\mbox{}+\pi H_i(X_{R+1},\ldots,X_8)
\end{eqnarray}
for $i=1,2$, with appropriate quadratic forms $G_i,H_i$ and linear forms
$L_j^{(i)}$, all defined over $\ov$.  Then if $\overline{H_1}$ and
$\overline{H_2}$ have a common zero over $F_v$ the pair $Q_1,Q_2$ is
not minimized.
\end{lemma}

To prove this we make a change of variables among $X_{R+1},\ldots,X_8$
so as to suppose that
\[\overline{H_1}(0,\ldots,0,1)=\overline{H_2}(0,\ldots,0,1)=0.\]
One then sets $T=\diag(\pi,\ldots,\pi,1)$ and
$U=\diag(\pi^{-2},\pi^{-2})$.  These satisfy (\ref{mincon}) and
produce a pair of integral forms $(Q_1,Q_2)_T^U$.  Thus $Q_1$ and
$Q_2$ cannot be minimized.

\section{The Case $r\le 4$}

Lemma \ref{notmin} shows that a minimized pair of forms cannot have
$R(q_1,q_2)\le 3$.  Our proof of Theorem \ref{T2} proceeds via a
case-by-case analysis of the remaining possibilities for $R(q_1,q_2)$
and $r(q_1,q_2)$.  In this section we examine the possibility that
$r(q_1,q_2)\le 4$.  We begin by proving the following result. 

\begin{lemma}\label{R4}
Suppose that $Q_1,Q_2$ is a minimized
pair of quadratic forms in 8 variables, with a common non-trivial zero
over $k_v$.  Then $R(q_1,q_2)\not=4$.
\end{lemma}

For the proof we argue by contradiction.
If $R(q_1,q_2)=4$ there is a change of variables in $\GL_8(\ov)$ so
that $q_1$ and $q_2$ are functions of $X_1,\ldots,X_4$ only.  Thus we
may write $Q_1$ and $Q_2$ in the shape (\ref{x8shape}) with $R=4$. In
view of Lemma \ref{x8} we know that
$\overline{H_1}$ and $\overline{H_2}$ cannot have a common zero over $F_v$.
We now claim similarly that the reductions
$\overline{G_1}$ and $\overline{G_2}$ cannot have a common zero over $F_v$.
If they did, then after a change of variables in $\GL_8(\ov)$ we could 
assume that
\[\overline{G_1}(1,0,0,0)=\overline{G_2}(1,0,0,0)=0.\]
Taking $T=\diag(1,\pi,\pi,\pi,1,1,1,1)$ it would then follow that both
the forms $Q_1(T\b{X})$ and $Q_2(T\b{X})$ are divisible by $\pi$.
Hence if $U=\diag(\pi^{-1},\pi^{-1})$ we would find that $(Q_1,Q_2)_T^U$ is a
pair of integral forms.  However the condition (\ref{mincon}) is
satisfied, since $\det(T)=|\pi|_v^3$ and $\det(U)=|\pi|_v^{-2}$. Since
$Q_1$ and $Q_2$ are minimized this contradiction proves our claim.

We therefore conclude that if the pair $Q_1,Q_2$ is minimized then neither
$\overline{G_1}=\overline{G_2}=0$ nor $\overline{H_1}=\overline{H_2}=0$
can have a non-trivial solution over $F_v$.  However it is then easy
to see that $Q_1=Q_2=0$ cannot have a non-trivial zero over $k_v$,
giving us the desired contradiction.  This completes our proof of
Lemma \ref{R4}.
\bigskip

We are now ready to show that the case $r(q_1,q_2)\le 4$ cannot occur.

\begin{lemma}\label{rge5}
Suppose that $\# F_v\ge 8$ and that  
$Q_1,Q_2$ is a minimized
pair of quadratic forms in 8 variables, with a common non-trivial zero
over $k_v$.  Then $r(q_1,q_2)\ge 5$.
\end{lemma}

We argue by contradiction.  We know from Lemmas \ref{notmin} and
\ref{R4} that $R=R(q_1,q_2)\ge 5$.  Thus if $r=r(q_1,q_2)\le 4$ we see
that either $r<R$ with $r$ even, or $r\le R-2$.  We may therefore
deduce from Lemma \ref{r-2} that
\[q_1=q_3(X_1,\ldots,X_{R-2})+X_{R-1}\ell(X_1,\ldots,X_{R-1})\]
and
\[q_2=q_4(X_1,\ldots,X_{R-2})+X_{R-1}X_R\]
with $r(q_3,q_4)=r-2\le 2$ and $R(q_3,q_4)\le R-2$.  
If $R(q_3,q_4)\le 2$ we may
make a change of variable so that 
\[q_1=q_3'(X_1,X_{2})+X_{R-1}\ell(X_1,\ldots,X_{R-1})\]
and
\[q_2=q_4'(X_1,X_{2})+X_{R-1}X_R.\]
These are in the form (\ref{3L}) and hence the pair $Q_1,Q_2$
cannot be minimized, by Lemma
\ref{notmin}.  This contradiction shows that $R(q_3,q_4)\ge 3$.

We may now repeat the previous argument, but applied to the forms
$q_3$ and $q_4$.  Since $r'=r(q_3,q_4)\le 2$ and $R'=R(q_3,q_4)\ge 3$
we must either have $r'<R'$ with $r'$ even, or $r'\le R'-2$.  In
either case we deduce from Lemma \ref{r-2} firstly that $q_4$ can 
be put into the shape
(\ref{q2sh}), then that $r'=2$ and finally that
$q_3=X_{R'-1}\ell'(X_1,\ldots,X_{R'-1})$ and $q_4=X_{R'-1}X_R'$.
It follows that $q_1$ and $q_2$ can be put into the form (\ref{3L})
(indeed with two terms on the right, rather than three).  As before,
Lemma \ref{notmin} implies that $Q_1,Q_2$ cannot be minimized, giving
the required contradiction.  This proves the lemma.

\section{The Case $R=5$}

So far we have shown that certain small values of $R(q_1,q_2)$ or
$r(q_1,q_2)$ are impossible. For large values we will show that the
pencil $<Q_1,Q_2>$ contains a form which splits off three
hyperbolic planes.  We have already remarked in Lemma \ref{r7}
that this is the case when
$r(q_1,q_2)\ge 7$.  In this section we deal with the case $R=5$.

\begin{lemma}\label{R5}
Suppose that $\# F_v\ge 9$.  Then if $R(q_1,q_2)=5$ there is at least
one form in the pencil 
$<Q_1,Q_2>$ which splits off three hyperbolic planes.
\end{lemma}

We begin by writing our forms in the shape (\ref{x8shape}) with
$R=5$. We know from Lemma \ref{rge5} 
that if $R=5$ we will also have $r=5$, so that if $g_i=\overline{G_i}$
for $i=1,2$, then $r(g_1,g_2)=5$.  

We claim that if $h_i=\overline{H_i}$ for $i=1,2$, then $R(h_1,h_2)=
r(h_1,h_2)=3$.  To prove this, suppose firstly that $R(h_1,h_2)\le
2$.  Then $h_1,h_2$ would have a common zero over $F_v$ and Lemma
\ref{x8} would contradict the minimality of the pair $Q_1,Q_2$.  
On the other hand, if $R(h_1,h_2)=3$ and $r=r(h_1,h_2)\le 2$ then
either $r$ is even or $r\le R(h_1,h_2)-2$.  We may therefore apply
Lemma \ref{r-2}.  Here the forms corresponding
to $q_3$ and $q_4$ must vanish identically since $r\le 2$.  It follows
that $h_1$ and $h_2$ have a common zero at $(0,0,1)$, after the change
of variables given in Lemma \ref{r-2}.  As above this leads, via Lemma
\ref{x8}, to a contradiction.  Thus we must have $R(h_1,h_2)=
r(h_1,h_2)=3$, as claimed.

Since $r(g_1,g_2)=R(g_1,g_2)=5$, the form $g_1+tg_2$ cannot be
singular for all values of $t$. However its determinant (or half-determinant) 
is a polynomial of degree at most 5, and we therefore see that 
$g_1+tg_2$ is singular for at most 5 values of $t$.  Similarly we find that
$h_1+th_2$ is singular for at most 3 values of $t$.  Since $\#
F_v>8$ there is a linear combination for which $g_1+tg_2$ has rank 5
and $h_1+th_2$ has rank 3. Let $\tau$ be any lift of $t$ to $\ov$, and
let $Q=Q_1+\tau Q_2$ so that
\[Q(X_1,\ldots,X_8)=G(X_1,\ldots,X_5)+\pi\sum_{i=1}^5X_iL_i(X_6,X_7,X_8)
+\pi H(X_6,X_7,X_8)\]
with $\rank(g)=5$ and $\rank(h)=3$.  We may then write
\[g(X_1,\ldots,X_5)=X_1X_2+X_3X_4+cX_5^2\;\;\;\mbox{and}\;\;\;
h(X_6,X_7,X_8)=X_6X_7+c'X_8^2\]
after a suitable change of variables.

Thus $Q$ is in the appropriate shape to apply Lemma \ref{lift} with $s=2$ and
\[\widetilde{Q}(X_5,\ldots,X_8)=cX_5^2+\pi X_5L(X_6,X_7,X_8)+\pi
H(X_6,X_7,X_8)\]
for a suitable linear form $L(X_6,X_7,X_8)$. We see from Lemma
\ref{lift} that $Q$ splits off two hyperbolic planes and so it remains to
show that $Q_0$ splits off at least one hyperbolic plane.  Let
$J(X_5,\ldots,X_8)=\pi^{-1}Q_0(\pi X_5,X_6,X_7,X_8)$.  Then $J$ has
coefficients in $\ov$ and satisfies $J\equiv H\modd{\pi}$.  Now, since 
$h(X_6,X_7,X_8)=X_6X_7+c'X_8^2$ we are able to make a second
application of Lemma \ref{lift} to show that $J$ splits off a
hyperbolic plane.  Lemma \ref{R5} then follows.

\section{The Case $R=8$}

In this section we prove the following result.

\begin{lemma}\label{R8}
Suppose that $\# F_v\ge 9$.  
Let $Q_1,Q_2$ be a minimized pair of forms in 8 variables, with a 
common zero over $k_v$. Then if $R(q_1,q_2)=8$ there is at least
one form in the pencil 
$<Q_1,Q_2>$ which splits off three hyperbolic planes.
\end{lemma}

The result is an immediate consequence of Lemma \ref{r7} unless
$r(q_1,q_2)\le 6$ as we henceforth suppose. It will be convenient to
set $r=r(q_1,q_2)$.  Then Lemma \ref{r-2} will apply, giving us
representations
\[q_1=q_3(X_1,\ldots,X_6)+X_7\ell(X_1,\ldots,X_6)\]
and
\[q_2=q_4(X_1,\ldots,X_6)+X_7X_8\]
with $\rank(q_3)=r(q_3,q_4)=r-2\le 4$ and
$5\le R(q_3,q_4)\le 6$.  Hence we may apply Lemma 
\ref{r-2} a second time, leading to expressions
\[q_1=q_5(X_1,\ldots,X_4)+X_5\ell'(X_1,\ldots,X_4)+X_7\ell(X_1,\ldots,X_6)\]
and
\[q_2=q_6(X_1,\ldots,X_4)+X_5X_6+X_7X_8\]
with $\rank(q_5)=r(q_5,q_6)\le 2$.  If $R(q_5,q_6)\ge 3$ we would be
able to repeat the process a third time, showing that $r(q_5,q_6)=2$
and producing
\[q_1=X_3\ell''(X_1,X_2)+X_5\ell'(X_1,\ldots,X_4)+
X_7\ell(X_1,\ldots,X_6)\]
and
\[q_2=X_3X_4+X_5X_6+X_7X_8.\]
However the forms $q_1,q_2$ would then be in the shape (\ref{3L}), so
that, according to Lemma \ref{notmin}, the pair $Q_1,Q_2$ could 
not have been minimized.  This contradiction shows that $R(q_5,q_6)\le
2$, allowing us to write
\[q_1=q_5(X_1,X_2)+X_5\ell'(X_1,\ldots,X_4)+X_7\ell(X_1,\ldots,X_6)\]
and
\[q_2=q_6(X_1,X_2)+X_5X_6+X_7X_8,\]
after a suitable change of variables.  It is now apparent that, in
order to have $R(q_1,q_2)=8$, the linear forms 
\[X_1,\;X_2,\;X_5,\;X_6,\;X_7,\;X_8,\;\ell'(X_1,\ldots,X_4)\;\;\mbox{and}\;\;
\ell(X_1,\ldots,X_6)\]
must be linearly independent.  We may therefore write
\beql{R8e}
\begin{array}{ccc}
q_1&=&q_5(X_1,X_2)+X_5X_3+X_7X_4\\
q_2&=&q_6(X_1,X_2)+X_5X_6+X_7X_8, \end{array}
\eeq
after a further change of variables.
\bigskip

To complete the argument we now call on the following result, which will
be used repeatedly in the rest of the proof of Theorem \ref{T2}.
\begin{lemma}\label{thelemma}
Let $Q_1(X_1,\ldots,X_8)$ and $Q_2(X_1,\ldots,X_8)$ be forms
over $\ov$, not necessarily minimized.  Suppose that 
\[Q_1(X_1,X_2,X_3,X_4,X_5,0,0,0)\equiv X_1X_2+X_3X_4\modd{\pi}\]
and that 
\[\pi\nmid Q_2(0,0,0,0,1,0,0,0).\]
Then there is at least one form in the pencil 
$<Q_1,Q_2>$ which splits off three hyperbolic planes.
\end{lemma}

We shall prove the lemma in a moment, but first we demonstrate how it may
be used to complete our proof of Lemma \ref{R8}.  It is not possible
for both the forms $q_5$ and $q_6$ in (\ref{R8e}) to vanish
identically, since $R(q_1,q_2)=8$.  Moreover, if both $q_5$ and $q_6$
were merely multiples of $X_1X_2$ at least one of $q_1$ or $q_2$ would
be a sum of three hyperbolic planes, in which case an application of Lemma 
\ref{lift} completes the proof.  We may therefore assume that $q_6$,
say, contains a non-zero term in $X_1^2$.  We now replace $Q_1$ by
$Q_1+cQ_2$ for a suitable $c\in\ov$ and substitute $X_3$ and $X_4$ for
$X_3+cX_6$ and $X_4+cX_8$ respectively.  This enables us to assume
that $q_5(1,0)=0$. Then, re-labelling the
variables, we may write
\[\begin{array}{ccc}
q_1&=&X_1X_2+X_3X_4+q_5(X_5,X_6)\\
q_2&=&X_1X_7+X_3X_8+q_6(X_5,X_6),\end{array}\]
with $q_5(1,0)=0$ and $q_6(1,0)\not=0$.  Thus Lemma \ref{thelemma} 
applies to $Q_1,Q_2$, and completes the proof of Lemma \ref{R8}
\bigskip

It therefore remains to establish Lemma \ref{thelemma}.  It will be
convenient to write 
\[S_i(X_1,X_2,X_3,X_4,X_5)=Q_i(X_1,X_2,X_3,X_4,X_5,0,0,0),\;\;\;(i=1,2).\]
We claim that there exists $\lambda\in \ov$ and $T\in\GL_5(\ov)$ 
such that 
\[T(0,0,0,0,1)=(0,0,0,0,1)\]
and
\beql{R8i}
(S_1-\lambda S_2)(T\b{X})=X_1X_2+X_3X_4.
\eeq
In particular $S_1-\lambda S_2$ will be singular.  Since $\ov$ and
$\GL_5(\ov)$ are compact it will suffice to show that for every
positive integer $f$ there are suitable $\lambda$ and $T$ such that
\beql{fc}
(S_1-\lambda S_2)(T\b{X})\equiv X_1X_2+X_3X_4\modd{\pi^f}.
\eeq
We will prove this by induction on $f$, the case $f=1$ being handled
by the hypotheses of the lemma.

We therefore suppose that (\ref{R8i}) holds for some particular $f$
and show how to derive a corresponding statement with exponent
$f+1$.  It will be convenient to write
\begin{eqnarray*}
S(\b{X})&=&(S_1-\lambda S_2)(T\b{X})\\
&=&X_1X_2+X_3X_4+\pi^f S'(X_1,\ldots,X_4)\\
&&\hspace{2cm}\mbox{}
+\pi^f L(X_1,\ldots,X_4)X_5+\pi^fcX_5^2
\end{eqnarray*}
and 
\[U(\b{X})=S_2(T\b{X})=U_0(X_1,X_2,X_3,X_4)+M(X_1,\ldots,X_4)X_5+dX_5^2,\]
Since $T(0,0,0,0,1)=(0,0,0,0,1)$ we have
$U(0,0,0,0,1)=S_2(0,0,0,0,1)$, whence $\pi\nmid d$.

We now examine $S-\pi^fcd^{-1}U$, which will have
coefficients in $\ov$.  By construction this form will have no term in
$X_5^2$.  Let
\begin{eqnarray*}
\lefteqn{(S-\pi^fcd^{-1}U)(X_1,X_2,X_3,X_4,0)=V(X_1,X_2,X_3,X_4)}\\
&=&
X_1X_2+X_3X_4+\pi^f
S'(X_1,\ldots,X_4)-\pi^fcd^{-1}U(X_1,X_2,X_3,X_4,0),
\end{eqnarray*}
say, so that
$V(X_1,X_2,X_3,X_4)\equiv X_1X_2+X_3X_4\modd{\pi}$. We therefore see that Lemma
\ref{lift} applies, producing a transform in $T_0\in\GL_4(\ov)$ such that
$V(T_0\b{X})=X_1X_2+X_3X_4$. Thus there is an admissible $T\in
GL_5(\ov)$ with
\[(S-\pi^fcd^{-1}U)(T\b{X})=X_1X_2+X_3X_4+\pi^fL'(X_1,\ldots,X_4)X_5,\]
where we have set
\[L'(X_1,\ldots,X_4)=L\big(T_0(X_1,\ldots,X_4)\big)-
cd^{-1}M\big(T_0(X_1,\ldots,X_4)\big).\]
We now make a
further change of variable, of the shape
\[X_i\rightarrow X_i+\pi^f\gamma_i X_5\;\;\;(1\le i\le 4)\]
to put $S-\pi^fcd^{-1}U$ into the form
\[(S-\pi^fcd^{-1}U)(T\b{X})=X_1X_2+X_3X_4+\pi^{2f}c'X_5^2.\]
Since $2f>f$ this establishes (\ref{fc}) with exponent $f+1$, and with
a new value $\lambda+\pi^fcd^{-1}$ in place of $\lambda$.  We have
therefore completed the induction step, thereby establishing the claim that
we can choose $\lambda$ so that (\ref{R8i}) holds.
\bigskip

Returning to the statement of Lemma \ref{thelemma} we now see that
there is a change of variable putting $Q_1-\lambda Q_2$ into the shape
\[X_1X_2+X_3X_4+X_6L_1(X_1,\ldots,X_8)+X_7L_2(X_1,\ldots,X_8)
+X_8L_3(X_1,\ldots,X_8).\]
If $\rank(Q_1-\lambda Q_2)=8$ then $X_5$ must appear in at least one 
of the forms $L_i$.  In this case after a further change of variable we 
can represent $Q_1-\lambda Q_2$ as
\[X_1X_2+X_3X_4+X_5X_6+X_7L_2'(X_1,\ldots,X_8)
+X_8L_3'(X_1,\ldots,X_8).\]
Another change of variable, of the form
\[X_i\rightarrow X_i+\mu_iX_7+\nu_iX_8,\;\;\; (1\le i\le 6)\]
then produces
\[Q_1-\lambda Q_2=X_1X_2+X_3X_4+X_5X_6+X_7L_2''(X_7,X_8)
+X_8L_3''(X_7,X_8),\]
giving us a form which splits off three hyperbolic planes.

It remains to consider the possibility that $Q_1-\lambda Q_2$ is
singular.  Here we use an argument shown to the author by Colliot-Th\'el\`ene.
We recall from Lemma \ref{LBP} that the form
$F(x,y)=\det(xM(Q_1)+yM(Q_2))$ cannot vanish, and will have distinct factors.
However 
\[F(1,-\lambda)=\det(Q_1-\lambda Q_2)=0,\]
whence
\[F(x,y)=(\lambda x+y)G(x,y)\]
for some form $G(x,y)\in k_v[x,y]$ with $G(1,-\lambda)\not=0$.
Suppose that we have $v(G(1,-\lambda))=e$, so that $\pi^{-e}G(1,-\lambda)$ is
a unit in $\ov$. Then if $r$ is a large enough integer we will have 
$v(G(1,\pi^r-\lambda))=e$.  Moreover
\[F(1,\pi^r-\lambda)=\pi^rG(1,\pi^r-\lambda),\]
whence
$v(F(1,\pi^r-\lambda))$ will be odd for any sufficiently large integer
$r$ of opposite parity to $e$.  In particular $F(1,\pi^r-\lambda)$ is
not a square in $k_v$ for such a choice of $r$, so that 
$Q_1+(\pi^r-\lambda)Q_2$ is nonsingular
and contains three hyperbolic planes, by virtue of Lemma \ref{nsq}.  This
completes the proof of Lemma \ref{thelemma}.
\bigskip

We end this section by giving a useful corollary to Lemma
\ref{thelemma}.
\begin{lemma}\label{tl+}
Let $Q_1(X_1,\ldots,X_8)$ and $Q_2(X_1,\ldots,X_8)$ be forms
over $\ov$, not necessarily minimized.  Suppose that 
\[Q_1(X_1,X_2,X_3,X_4,X_5,0,0,0)\equiv X_1X_2+X_3X_4\modd{\pi}\]
and
\[Q_2(X_1,X_2,X_3,X_4,X_5,0,0,0)\equiv Q_2'(X_1,X_2,X_3,X_4)\modd{\pi}\]
for some quaternary form $Q_2'$.  Suppose further that 
\[\pi^2\mid Q_1(0,0,0,0,1,0,0,0)\]
and
\[\pi^2\nmid Q_2(0,0,0,0,1,0,0,0)\]
Then there is at least one form in the pencil 
$<Q_1,Q_2>$ which splits off three hyperbolic planes.
\end{lemma}
If $\pi\nmid Q_2(0,0,0,0,1,0,0,0)$ this follows at once from Lemma
\ref{thelemma}. Otherwise we define $(Q_1,Q_2)_T^U=(V_1,V_2)$ with
$T=\diag(\pi,\pi,\pi,\pi,1,\pi^3,\pi^3,\pi^3)$ and
$U=\diag(\pi^{-2},\pi^{-1})$. Then $V_1,V_2$ are integral forms.
Moreover we have
\[V_1(\b{X})\equiv X_1X_2+X_3X_4+L(X_1,X_2,X_3,X_4)X_5+cX_5^2\modd{\pi}\]
for some linear form $L$ and some $c\in\ov$, while
\[V_2(\b{X})\equiv dX_5^2\modd{\pi}\]
for some unit $d\in\ov$.  A suitable substitution $X_i\rightarrow
X_i+c_iX_5$ for $1\le i\le 4$ then transforms $V_1$ so that
\[V_1(\b{X})\equiv X_1X_2+X_3X_4+c'X_5^2\modd{\pi}\]
say, while leaving $V_2(\b{X})\equiv dX_5^2\modd{\pi}$.  We may then
apply Lemma \ref{thelemma} to the forms $V_1-c'd^{-1}V_2$ and $V_2$,
and Lemma \ref{tl+} follows.

\section{The Case $\rank(q_2)\le 2$ --- First Steps}

In the next two sections we examine the case in which $q_2$ has
small rank.  This turns out to be an important prelude to our
treatment of pairs for which $R(q_1,q_2)=6$.  Our goal is the
following result.

\begin{lemma}\label{rk2f}
Suppose that $\# F_v\ge 9$.  
Let $Q_1,Q_2$ be a minimized pair of forms with a common zero over
$k_v$. Then if $\rank(q_2)\le 2$ there is at least
one form in the pencil 
$<Q_1,Q_2>$ which splits off three hyperbolic planes.
\end{lemma}

However in the present section we will content ourselves with the
following intermediate statement.

\begin{lemma}\label{rk2}
Suppose that $\# F_v\ge 9$.  
Let $Q_1,Q_2$ be a minimized pair of forms over
$k_v$ such that $\rank(q_2)\le 2$.  Suppose indeed that
$q_2(\b{X})=q_2(X_1,X_2)$. Then either there is at least
one form in the pencil 
$<Q_1,Q_2>$ which splits off three hyperbolic planes, or $q_2$ is an
anisotropic form of rank 2, and there is a linear
change of variable in $\GL_6(\ov)$, involving only $X_3,\ldots,X_8$
which makes 
\beql{falt}
q_1(\b{X})=X_1X_5+X_2X_6+n(X_3,X_4).
\eeq
\end{lemma}

We begin by observing that we must have $\rank(q_2)=2$.  Indeed we
will be able to write $q_2(\b{X})=n(X_1,X_2)$ for some anisotropic form $n$,
after a suitable change of variables. To see this we notice that in
all other cases we can write $q_2$ as a product of linear factors over $F_v$.
Thus, after a change of variable we can take
$q_2(\b{X})=X_1\ell(\b{X})$.  The transforms
$T=\diag(\pi,1,\ldots,1)$ and $U=(\pi^{-1},1)$ then make
$(Q_1,Q_2)_T^U$ integral, but would satisfy (\ref{mincon}), contradicting
the minimality of the pair $Q_1,Q_2$.

In view of Lemma \ref{r7} we can assume that $r=r(q_1,q_2)\le
6$. Similarly, in view of Lemmas \ref{notmin}, \ref{R4} and \ref{R5}
we may suppose that $R=R(q_1,q_2)\ge 6$.  We
proceed to use these conditions to narrow down the possible shapes that
$q_1$ may take.  With $q_2(\b{X})=n(X_1,X_2)$, we write 
\[q_1(\b{X})=q_3(X_1,X_2)+X_1\ell_1(X_3,\ldots,X_8)+X_2\ell_2(X_3,\ldots,X_8)
+q_4(X_3,\ldots,X_8)\]
for appropriate quadratic forms $q_i$ and linear forms
$\ell_i$. Let $\rank(q_4)=m$, say, and change variables so as to write
\begin{eqnarray}\label{q1f}
q_1(\b{X})&=&q_3(X_1,X_2)+X_1\ell_1'(X_3,\ldots,X_{m+2})+
X_2\ell_2'(X_3,\ldots,X_{m+2})\nonumber\\
&&\hspace{5mm}\mbox{}
+X_1\ell_1''(X_{m+3},\ldots,X_8)+X_2\ell_2''(X_{m+3},\ldots,X_8)\nonumber\\
&&\hspace{1cm}\mbox{}+q_5(X_3,\ldots,X_{m+2})
\end{eqnarray}
with $\rank(q_5)=m$.
Our analysis now splits into 3 cases in which the linear forms
$\ell_1'',\ell_2''$ both vanish, or are
linearly dependent but do not both vanish, or are linearly independent.

\subsection{Case 1.}  
If $\ell_1''$ and $\ell_2''$ both vanish identically then
$R(q_1,q_2)\le m+2$. Since we are assuming that $R\ge 6$ this shows
that $m\ge 4$.  On the other hand, if $\chi(F_v)\not=2$, or
$m$ is even, the
determinant of $q_1+tq_2$ (considered as a quadratic form in $m+2$
variables) is a polynomial in $t$ of degree at most $m+2$, in which the term
in $t^2$ has coefficient $\det(n)\det(q_5)\not=0$.  It follows that
$r(q_1,q_2)\ge m+2$.  Similarly if $\chi(F_v)=2$ and $m$ is odd, 
the half-determinant of 
$q_1+tq_2$ is a polynomial in which the coefficient of $t^2$ is 
$\det(n)\hd(q_5)\not=0$,
and again we conclude that $r(q_1,q_2)\ge m+2$.  We are supposing that
$r\le 6$, so we must have $m\le 4$.

We are therefore left with the case in which $m=4$, in which situation
a change of variables will allow us to put $q_5$ into one of the shapes
$X_3X_4+X_5X_6$ or $X_3X_4+n(X_5,X_6)$. We can now use further linear
transformations, of the type 
\[X_i\rightarrow X_i+\lambda_i(X_1,X_2)\;\;\;(3\le i\le 6)\]
to put $q_1$ into one of the shapes
\beql{2Nc}
q_1(\b{X})=q_6(X_1,X_2)+X_3X_4+X_5X_6
\eeq
or
\beql{2Nc'}
q_1(\b{X})=q_6(X_1,X_2)+X_3X_4+n(X_5,X_6).
\eeq
\bigskip

We now call on the following lemma, which we shall prove shortly.
\begin{lemma}\label{2N}
Let $s_1(X,Y)$ and $s_2(X,Y)$ be quadratic forms over a finite field
$F$. Suppose that $s_1$ and $s_2$ have no common factor and that 
$r(s_1,s_2)=2$.  Then there are at least $\tfrac12 (\# F-1)^2$ 
pairs $a,b\in F$, not both zero, for which
$as_1+bs_2$ is a hyperbolic plane, and at least $\tfrac12 (\# F-1)^2$ 
such pairs for which it is anisotropic of rank 2.
\end{lemma}
\bigskip

It is clear that $r(q_6,n)=2$, because $\rank(n)=2$.  Moreover, since $n$
cannot have a linear factor, the forms $q_6$ and $n$ are coprime unless
$q_6$ is a multiple of $n$.  
Suppose firstly that $q_6$ is not a multiple of $n$.  In
Lemma \ref{2N} we must have $\tfrac12
(\# F_v-1)^2>\# F_v-1$ if $\# F_v\ge 9$, and we therefore conclude that 
there are linear
combinations $aq_6+bn$ with $a\not=0$ which are hyperbolic planes, and also
which are anisotropic of rank 2. We may then deduce that 
there is 
a linear combination $q_1'=q_1+cq_2$ taking one of the forms
\[X_1X_2+X_3X_4+X_5X_6,\;\;\;  \mbox{or}\;\;\; n(X_1,X_2)+X_3X_4+n(X_5,X_6),\]
in the two cases given by (\ref{2Nc}) and (\ref{2Nc'}).  As explained in Section
\ref{1QF}, if $q_1'(\b{X})=n(X_1,X_2)+X_3X_4+n(X_5,X_6)$ then a change
of variable will make $q_1'$ a sum of 3 hyperbolic planes, and so in either case
Lemma \ref{rk2} follows from Lemma \ref{lift}.

In the alternative situation in which $q_6$ is a multiple of $n$, there is 
a linear combination $q_1'=q_1+cq_2$ taking one of the forms
\[X_3X_4+X_5X_6,\;\;\;  \mbox{or}\;\;\; n(X_1,X_2)+X_3X_4+n(X_5,X_6),\]
in the two cases given by (\ref{2Nc}) and (\ref{2Nc'}) respectively.  
As above, in the second of
these cases $q_1'$ is a sum of 3 hyperbolic planes, and again we may
complete the proof of Lemma \ref{rk2}.  
We may therefore reduce our considerations to the case in which
\[q_1(\b{X})=X_3X_4+X_5X_6,\;\;\;  \mbox{and}\;\;\;  q_2(\b{X})=n(X_1,X_2).\]
However in this situation we may apply Lemma \ref{thelemma} 
after re-ordering the variables, since $\pi\nmid n(1,0)$.  This
completes the treatment of Case 1.
\bigskip

We now present our proof of Lemma \ref{2N}.  
It will be convenient to write $\# F=N$ temporarily.
Let $\cl{S}$ be the set of quadruples $(a,b,x,y)\in F$ for which
$(a,b)\not=(0,0)$ and $(x,y)\not=(0,0)$ and such that $as_1(x,y)+bs_2(x,y)=0$.
Since $s_1$ and $s_2$ have no factor in common there is no
$(x,y)\not=(0,0)$ such that $s_1(x,y)=s_2(x,y)=0$.  Thus there are $N-1$ 
pairs $(a,b)$ for each admissible pair $(x,y)$.  It follows that
\[\#\cl{S}=(N-1)(N^2-1).\]

On the other hand, if $(a,b)$ is a pair such that
$as_1(X,Y)+bs_2(X,Y)$ is a hyperbolic plane
then there are $2(N -1)$ corresponding pairs $(x,y)$, while if 
$as_1(X,Y)+bs_2(X,Y)$ is anisotropic of rank 2 there are
none. Moreover, in the remaining case, in which $as_1(X,Y)+bs_2(X,Y)$
has a repeated linear factor, there are $N-1$ pairs $x,y$.  Suppose
that the hyperbolic plane, anisotropic, and repeated factor cases occur
$N_h, N_a$ and $N_r$ times each respectively.  Then 
\beql{2Ni}
N_h+N_a+N_r=N^2-1
\eeq
and
\beql{2Nii}
2(N-1)N_h+(N-1)N_r=\#\cl{S}=(N-1)(N^2-1).
\eeq
However the determinant of $as_1(X,Y)+bs_2(X,Y)$ is a binary quadratic
form in $a$ and $b$ which does not vanish identically, since
$r(s_1,s_2)=2$. It follows that $N_r\le 2(N-1)$.  We therefore
deduce from (\ref{2Ni}) and (\ref{2Nii}) that
\[2N_h=N^2-1-N_r\ge N^2-1-2(N-1)\]
and that
\[2N_a=2(N^2-1-N_h-N_r)=N^2-1-N_r\ge N^2-1-2(N-1).\]
These inequalities suffice for the lemma.

\subsection{Case 2.}  
The next case is that in which $\ell_1''$ and $\ell_2''$ are not both zero
but are linearly dependent.  After a change of variable between $X_1$
and $X_2$ we may suppose that $\ell_2''=0$, and then, after a further
change of variable among $X_{m+3},\ldots,X_8$ we may assume that
$\ell_1''=X_{m+3}$.  It follows that $R(q_1,q_2)\le m+3$.  Since we are
assuming that $R\ge 6$ this shows that $m\ge 3$.
The representation (\ref{q1f}) then becomes
\begin{eqnarray}\label{C2q1}
q_1(\b{X})&=&
q_4(X_1,X_2)+X_1\ell_3(X_3,\ldots,X_{m+2})+X_2\ell_4(X_3,\ldots,X_{m+2})\nonumber\\
&&\hspace{2cm}\mbox{}+X_1X_{m+3}+q_5(X_3,\ldots,X_{m+2}).
\end{eqnarray}
(Our numbering of forms $q_i$ and $\ell_i$ will be independent of that
used for Case~1.)  
Thus if the coefficients of $X_2^2$ in $q_1$ and $q_2$ are $a$ and
$b\not=0$ respectively, we may write $q_1+tq_2$ as
\beql{1t2}
q_1+tq_2=q_6(X_2,\ldots,X_{m+2})+X_1(X_{m+3}+\ell(X_1,\ldots,X_{m+2};t))
\eeq
for some linear form $\ell(X_1,\ldots,X_{m+2};t)$ depending on $t$,
where
\[q_6=(a+tb)X_2^2+X_2\ell_4(X_3,\ldots,X_{m+2})+q_5(X_3,\ldots,X_{m+2}).\]
Then $\rank(q_6)\ge \rank(q_5)=m$ and so $r(q_1,q_2)\ge\rank(q_6)+2\ge
m+2$, by (\ref{r+2}).  It follows that $m\le 4$.
In fact, if $m=4$ then the half-determinant of $q_6$
is a polynomial in $t$ whose linear term has a coefficient
$b\det(q_5)\not=0$.  Thus there is some value of $t$ for which $q_6$
has rank 5.  This however would imply that $q_1+tq_2$, given by
(\ref{1t2}), has rank 7, contradicting our assumption that
$r(q_1,q_2)\le 6$.  We must therefore have $m\le 3$, and indeed we can
deduce that $m=3$, since $6\le R(q_1,q_2)\le m+3$.

After a change of variable amongst $X_3,X_4,X_5$ we can now take
\[q_5(X_3,X_4,X_5)=X_3X_4+X_5^2,\]
whence (\ref{C2q1}) becomes
\begin{eqnarray*}
q_1(\b{X})&=&
q_4(X_1,X_2)+X_1\ell_5(X_3,X_4,X_5)+X_2\ell_6(X_3,X_4,X_5)\\
&&\hspace{1cm}\mbox{}+X_1X_6+X_3X_4+X_5^2.
\end{eqnarray*}
Further transforms of the type $X_3\rightarrow X_3+\lambda X_2$,
$X_4\rightarrow X_4+\mu X_2$ and $X_6\rightarrow X_6+\ell_7(X_1,\ldots,X_5)$
simplify this to
\[q_1(\b{X})=cX_2^2+dX_2X_5+X_1X_6+X_3X_4+X_5^2\]
for certain $c,d\in F_v$.  Now, with the same notation $b$ for the
coefficient of $X_2^2$ in $q_2=n(X_1,X_2)$ we consider
$q=q_1-cb^{-1}q_2$, which takes the shape
\[q=X_1(X_6+\ell_8(X_1,X_2))+dX_2X_5+X_3X_4+X_5^2.\]
After replacing $X_6+\ell_8(X_1,X_2)$ by $X_6$ we finally see that we can
replace our original forms $Q_1,Q_2$ by a pair $Q,Q_2$ for which,
after a change of variables in $\GL_8(\ov)$, we have
\[q=X_1X_6+dX_2X_5+X_3X_4+X_5^2\;\;\;\mbox{and}\;\;\;
q_2=n(X_1,X_2).\]
We may now re-label the variables so that
$q=X_1X_2+X_3X_4+dX_5X_6+X_6^2$ and $q_2=n(X_1,X_5)$.  Thus we can
apply Lemma \ref{thelemma} to complete the proof.

\subsection{Case 3.}  
If $\ell_1''$ and $\ell_2''$ are linearly independent we may take them to be 
$X_{m+3}$ and $X_{m+4}$ respectively, after an appropriate change of
variable. The representation (\ref{q1f}) then becomes
\begin{eqnarray*}
q_1(\b{X})&=&
q_3(X_1,X_2)+X_1\ell_1'(X_3,\ldots,X_{m+2})+X_2\ell_2'(X_3,\ldots,X_{m+2})\\
&&\hspace{1cm}\mbox{}+X_1X_{m+3}+X_2X_{m+4}+q_5(X_3,\ldots,X_{m+2}),
\end{eqnarray*}
which simplifies to
\[q_1(\b{X})=X_1X_{m+3}+X_2X_{m+4}+q_5(X_3,\ldots,X_{m+2})\]
after a substitution of the form 
\[X_{m+3}\rightarrow X_{m+3}+\ell_3(X_1,\ldots,X_{m+2}),\;\;\;
X_{m+4}\rightarrow X_{m+4}+\ell_4(X_1,\ldots,X_{m+2}).\]
It is clear
that $R(q_1,q_2)\le m+4$.  Since we are assuming that $R\ge 6$ this shows
that $m\ge 2$.  Thus $q_5$ will split off a hyperbolic plane unless $m=2$ and
$q_5$ is anisotropic. Thus either $q_1$ is a sum of three hyperbolic planes,
or we have
\[q_1(\b{X})=X_1X_5+X_2X_6+n(X_3,X_4).\]
Thus one or other of the conclusions in Lemma \ref{rk2} will hold.

\section{Completion of the Proof for $\rank(q_2)\le 2$}

We proceed with the proof of Lemma \ref{rk2f}.  Clearly we may assume
that $q_1$ takes the form (\ref{falt}) and that
$q_2(\b{X})=n(X_1,X_2)$.
We begin by setting
\beql{Sd}
S_i(X_7,X_8)=\pi^{-1}Q_i(0,0,0,0,0,0,X_7,X_8),\;\;\;(i=1,2).
\eeq
These forms are defined over $\ov$ since $q_1$ and $q_2$ contain no
terms in $X_7,X_8$.  Suppose firstly that $S_2$ does not vanish modulo
$\pi$.  In this case we can make a change of variables in $X_7$ and
$X_8$ so as to suppose that $\pi\nmid S_2(1,0)$. Then, with
$\lambda=S_1(1,0)S_2(1,0)^{-1}$, we replace $Q_1$ by
$Q=Q_1-\lambda Q_2$, which will contain no term in $X_7^2$.  We have
\[q(\b{X})=X_1X_5+X_2X_6+n(X_3,X_4)-\lambda n(X_1,X_2),\]
so that a change of variable, of the type
\[X_5\rightarrow X_5+\ell(X_1,X_2),\;\;\; X_6\rightarrow
X_6+\ell'(X_1,X_2)\]
puts $q(\b{X})$ into the shape $X_1X_5+X_2X_6+n(X_3,X_4)$ while
having no effect on $q_2(\b{X})=n(X_1,X_2)$.  These manoeuvres allow us
in effect to assume that $Q_1$ contains no term in $X_7^2$. If we now
relabel the variables so that $X_1,X_5,X_2,X_6,X_7$ become
$X_1,X_2,X_3,X_4,X_5$ respectively we see that Lemma \ref{tl+}
applies, showing that the pencil contains a form splitting off three
hyperbolic planes. 

We can therefore assume that $S_2$ vanishes modulo $\pi$.  Next, unless
$\overline{S_1}$ is anisotropic of rank 2 it will have a zero over
$F_v$, so that Lemma \ref{x8} may be applied.  Since this would
contradict the minimality of the pair $Q_1,Q_2$ we may therefore
assume that 
\[S_1(X_7,X_8)\equiv N(X_7,X_8),\;\; S_2(X_7,X_8)\equiv 0\modd{\pi}.\] 

We now consider the pair of forms $(Q_1',Q_2')=(Q_1,Q_2)_T^U$, where
\[T=\diag(\pi,\pi,1,\ldots,1)\;\;\;\mbox{and}\;\;\; 
U=\diag(1,\pi^{-1}).\]
Then $Q_1'$
and $Q_2'$ have coefficients in $\ov$, and $T$ and $U$ satisfy (\ref{TU=}).
Since the pair $Q_1,Q_2$ is minimized
we then deduce that $Q_1',Q_2'$ is also minimized.  However
$q_1'=n(X_3,X_4)$ so that Lemma \ref{rk2} applies.  We therefore conclude
that the pencil generated by $Q_1'$ and $Q_2'$ contains a form which
splits off 3 hyperbolic planes, except possibly when $q_2'$ 
takes the shape
\beql{q2'}
X_3\ell_1+X_4\ell_2+n(\ell_3,\ell_4)
\eeq
for linearly independent linear forms $\ell_i(X_1,X_2,X_5,X_6,X_7,X_8)$. 
However $q_2'$ contains no terms in $X_1$ or $X_2$, by construction.
Since $n$ is anisotropic we see that $\ell_3$ and $\ell_4$ must be independent
of $X_1$ and $X_2$, and then we deduce that $\ell_1$ and $\ell_2$ must also 
be independent of $X_1$ and $X_2$.  Similarly we observe that $q_2'$
contains no quadratic terms in $X_7$ and $X_8$, by construction,
whence $\ell_3$ and $\ell_4$ cannot involve $X_7$ or $X_8$.  Thus $\ell_1$ and
$\ell_2$ are functions of $X_5,\ldots,X_8$, while $\ell_3$ and $\ell_4$ are
functions of $X_5$ and $X_6$ alone.  Moreover $\ell_3$ and $\ell_4$ are
linearly independent.

We are finally in a position to complete the treatment of Lemma \ref{rk2f}.
The considerations above have shown that if $q_2'$ is given by
(\ref{q2'}) we may assume that our forms
take the shape
\begin{eqnarray*}
Q_1(\b{X})&\equiv& X_1X_5+X_2X_6+N(X_3,X_4)\\
&&\hspace{1cm}\mbox{}+\pi\{\widehat{Q_1}(X_1,\ldots,X_6)
+X_7M_1(X_1,\ldots,X_6)\\
&&\hspace{2cm}\mbox{}+X_8M_2(X_1,\ldots,X_6)+N(X_7,X_8)\}\modd{\pi^2}
\end{eqnarray*}
and
\begin{eqnarray*}
Q_2(\b{X})&\equiv & N(X_1,X_2)+
\pi\{X_3L_1(X_5,\ldots,X_8)+X_4L_2(X_5,\ldots,X_8)\\
&&\hspace{1cm}\mbox{}+X_1L_5(X_1,\ldots,X_8)+X_2L_6(X_1,\ldots,X_8)\\
&&\hspace{1cm}\mbox{}+N\big(L_3(X_5,X_6),L_4(X_5,X_6)\big)\}\;\modd{\pi^2},
\end{eqnarray*}
where $L_1,\ldots,L_4$ are lifts of $\ell_1,\ldots,\ell_4$.
According to our hypotheses there is a non-trivial point
$\b{x}\in k_v^8$ for which $Q_1=Q_2=0$.  By a suitable rescaling we
may suppose that $\b{x}$ has entries in $\ov$, at least one of which is
a unit. Since $\pi\mid Q_2(\b{x})$ we deduce that $\pi\mid
N(x_1,x_2)$, whence $\pi\mid x_1,x_2$.  Then, since $\pi\mid
Q_1(\b{x})$ we deduce that $\pi\mid N(x_3,x_4)$, leading to $\pi\mid
x_3,x_4$. We next use the fact that $\pi^2\mid Q_2(\b{x})$, which
shows that $\pi\mid N\big(L_3(X_5,X_6),L_4(X_5,X_6)\big)$.  This
implies $\pi\mid L_3(x_5,x_6),L_4(x_5,x_6)$, and since the forms $L_3$
and $L_4$ are independent modulo $\pi$ we find that $\pi\mid x_5,x_6$.
Finally, since $\pi^2\mid Q_1(\b{x})$ we obtain $\pi\mid N(x_7,x_8)$,
leading to $\pi\mid x_7,x_8$. We therefore have $\pi\mid\b{x}$,
contrary to our assumption.  Hence it is impossible for $q_2'$ to take
the shape (\ref{q2'}).  We therefore conclude that in every case the
pencil contains a form which splits off 3 hyperbolic planes.
This concludes the proof of Lemma \ref{rk2f}.

\section{The Case $R=6$}

We next suppose that $R(q_1,q_2)=6$, and prove the following result.
\begin{lemma}\label{R6}
Suppose that $\# F_v\ge 16$.  
Let $Q_1,Q_2$ be a minimized pair of forms with a common zero over
$k_v$. Then if $R(q_1,q_2)=6$ there is at least
one form in the pencil 
$<Q_1,Q_2>$ which splits off three hyperbolic planes.
\end{lemma}

In view of Lemma \ref{rge5} we may assume that $r(q_1,q_2)\ge 5$.
Since we have $R(q_1,q_2)=6$ we may write $q_i(\b{X})=q_i'(X_1,\ldots,X_6)$
for $i=1,2$.  We now consider the forms $S_i$ given by (\ref{Sd}). These will
have coefficients in $\ov$, so that we may define $s_i$ as the reduction
$\overline{S_i}$ over $F_v$.  The forms $s_1$ and $s_2$ cannot have a 
common zero over $F_v$ by Lemma \ref{x8}.

We first consider the possibility that
$s_1$ and $s_2$ are proportional. Then by considering appropriate
linear combinations of $Q_1$ and $Q_2$ we may suppose indeed that
$s_1=0$.  We now examine the forms $(V_1,V_2)=(Q_1,Q_2)_T^U$ where
$T=\diag(\pi,\ldots,\pi,1,1)$ and $U=\diag(\pi^{-2},\pi^{-1})$.  The
forms $V_1,V_2$ will be integral, while $T$ and $U$ satisfy
(\ref{TU=}).  Thus $V_1$ and $V_2$ are minimized.  Moreover, if $v_2$
is the reduction of $V_2$ to $F_v$ then $\rank(v_2)\le 2$. It then
follows from Lemma \ref{rk2f} that the pencil $<V_1,V_2>$ contains a
form which splits off three hyperbolic planes, and the required conclusion
then follows, subject to our assumption that $s_1$ and $s_2$ are
proportional. 

If $r(s_1,s_2)=1$, then every non-trivial
linear combination has rank at most one, and hence has a
zero. Otherwise, since we are now supposing that $s_1$ and $s_2$ are
not proportional and have no common zero, Lemma \ref{2N} shows that there are
at least $\tfrac12(\# F_v-1)^2$ non-trivial linear combinations $as_1+bs_2$
with a non-trivial zero.  In
contrast, at most $6(\# F_v-1)$ non-trivial linear combinations
$aq_1+bq_2$ can have rank smaller than $r(q_1,q_2)$, as one sees by
considering the corresponding determinant, or half determinant, as a
binary form in $a$ and $b$.  It follows that if $\# F_v>13$ then
there is at least one form $aq_1+bq_2$ with rank at least 5, such that
$as_1+bs_2$ has a non-trivial zero. By a change of variable between
$X_7$ and $X_8$ we may arrange that $as_1+bs_2$ vanishes at $(1,0)$.
Let us suppose, say, that $a\not=0$. Then $s_2(1,0)$ cannot vanish,
since $(1,0)$ cannot be a common zero of $s_1$ and $s_2$. Thus, on
replacing $Q_1$ by $aQ_1+bQ_2$ we may suppose that $\rank(q_1)\ge 5$,
that $s_1(1,0)=0$ and that $s_2(1,0)\not=0$. Since $\rank(q_1)\ge 5$
we know that $q_1$ splits off at
least two hyperbolic planes, whence we may make a change of 
variable so as to put
$q_1$ into the form $X_1X_2+X_3X_4+q_3(X_5,X_6)$. On interchanging the
variables $X_5$ and $X_7$ we now see that our forms are in the correct
shape for an application of Lemma \ref{tl+}, which provides a suitable
form splitting off 3 hyperbolic planes.

\section{Completion of the Proof of 
Theorem \ref{T2} --- The Case $R=7$}\label{secR7}

We turn now to the final case, in which $R(q_1,q_2)=7$.
\begin{lemma}\label{R7}
Suppose that $\# F_v\ge 32$.  
Let $Q_1,Q_2$ be a minimized pair of forms with a common zero over
$k_v$. Then if $R(q_1,q_2)=7$ there is at least
one form in the pencil 
$<Q_1,Q_2>$ which splits off three hyperbolic planes.
\end{lemma}

Our argument will depend on the following result, which we prove at the
end of this section.
\begin{lemma}\label{7l}
Let $q_1=Y_0Y_1+s_1(Y_1,Y_2,Y_3)$ and $q_2=Y_0Y_2+s_2(Y_1,Y_2,Y_3)$ be
quadratic forms over a finite field $F_v$ with $\# F_v\ge 32$. Suppose that
$q_1(0,0,1)$ and $q_2(0,0,1)$ are not both zero.  Then either
the forms have a singular common zero over $F_v$, 
or there are at least $5(\# F_v-1)$
non-zero pairs $(a,b)\in F_v^2$ for which $aq_1+bq_2$ is a
sum of two hyperbolic planes.
\end{lemma}
\bigskip

We now present our proof of Lemma \ref{R7}.
In view of Lemmas \ref{r7} and \ref{rge5} we may restrict attention to
the cases $r=5$ and $r=6$. Hence Lemma~\ref{r-2} shows that we may take 
\[q_1(\b{X})=q_3(X_1,\ldots,X_5)+X_6\ell(X_1,\ldots,X_6)\]
and
\[q_2(\b{X})=q_4(X_1,\ldots,X_5)+X_6X_7\]
with $\rank(q_3)=r(q_3,q_4)=r-2=3$ or $4$, and $R(q_3,q_4)=4$ or $5$.

We begin by considering the easy case in which $R(q_3,q_4)=5$.  Here we may
make a second application of Lemma \ref{r-2} and a further change 
of variables so as to write
\[q_1(\b{X})=q_5(X_1,X_2,X_3)+X_4\ell'(X_1,\ldots,X_4)+X_6\ell(X_1,\ldots,X_6)\]
and
\[q_2(\b{X})=q_6(X_1,X_2,X_3)+X_4X_5+X_6X_7\]
with $\rank(q_5)=r(q_5,q_6)=r-4=1$ or $2$. If $R(q_5,q_6)=3$ 
Lemma \ref{r-2} would show that
$q_5$ and $q_6$ have a common linear factor.  This however would allow
$q_1$ and $q_2$ to be written in the form (\ref{3L}), contradicting the
minimality of the pair $Q_1,Q_2$. Thus $R(q_5,q_6)\le 2$, so that we
can write
\[q_1(\b{X})=q_5(X_1,X_2)+X_4\ell'(X_1,\ldots,X_4)+X_6\ell(X_1,\ldots,X_6)\]
and
\[q_2(\b{X})=q_6(X_1,X_2)+X_4X_5+X_6X_7.\]
If we now set 
$T=\diag(\pi,\pi,1,\pi,1,\pi,1,1)$ and $U=\diag(\pi^{-1},\pi^{-1})$ we
see that $U$ and $T$ satisfy (\ref{TU=}), and that the forms
$(Q_1,Q_2)_T^U=(V_1,V_2)$, say, are integral.  Thus $V_1$ and $V_2$
are also minimized.  However one readily checks that if $v_1$ and
$v_2$ are the reductions of $V_1$ and $V_2$ to $F_v$ then
$R(v_1,v_2)\le 6$, since $v_1$ and $v_2$ contain no terms in $X_1$ and
$X_2$. Thus our previous results show that the pencil $<V_1,V_2>$
contains a form splitting off three hyperbolic planes, and this
suffices for the lemma.

We may therefore restrict our attention to the case in which
$R(q_3,q_4)=4$ and $r(q_3,q_4)=3$ or $4$.  Thus we may assume that
\[q_1(\b{X})=q_3(X_1,\ldots,X_4)+X_6\ell(X_1,\ldots,X_6)\]
and
\[q_2(\b{X})=q_4(X_1,\ldots,X_4)+X_6X_7.\]
Since $R(q_1,q_2)=7$ the variable $X_5$ must genuinely occur in
$\ell(X_1,\ldots,X_6)$ and so we can make a change of variable to obtain
\[q_1(\b{X})=q_3(X_1,\ldots,X_4)+X_6X_5\]
and
\[q_2(\b{X})=q_4(X_1,\ldots,X_4)+X_6X_7.\]

We proceed to write
\begin{eqnarray*}
Q_1(\b{X})&=&Q_3(X_1,\ldots,X_4)+X_5X_6\\
&&\hspace{1cm}\mbox{}+\pi\{Q_5(X_5,\ldots,X_8)+
\sum_{i=5}^8X_iL_i^{(1)}(X_1,\ldots,X_4)\}
\end{eqnarray*}
and similarly
\begin{eqnarray*}
Q_2(\b{X})&=&Q_4(X_1,\ldots,X_4)+X_6X_7\\
&&\hspace{1cm}\mbox{}+\pi\{Q_6(X_5,\ldots,X_8)+
\sum_{i=5}^8X_iL_i^{(2)}(X_1,\ldots,X_4)\},
\end{eqnarray*}
for suitable integral quadratic forms $Q_3,\ldots,Q_6$ and linear forms
$L_i^{(j)}$ such that the reductions $Q_3$ and $Q_4$ are $q_3$ and
$q_4$ respectively.
We note that $\pi$ cannot divide both $Q_5(0,0,0,1)$
and $Q_6(0,0,0,1)$, by Lemma \ref{x8}.  

We now examine $(G_1,G_2)=(Q_1,Q_2)_T^U$, where
\beql{TU1}
T=\diag(\pi,\pi,\pi,\pi,1,\pi,1,1)\;\;\;\mbox{and}\;\;\;
U=\diag(\pi^{-1},\pi^{-1}).
\eeq
Then 
\[G_1(\b{X})=\pi Q_3(X_1,\ldots,X_4)+H_1(X_5,\ldots,X_8)+
\pi\sum_{i=5}^8X_iM_i^{(1)}(\b{X})\]
and
\[G_2(\b{X})=\pi Q_4(X_1,\ldots,X_4)+H_2(X_5,\ldots,X_8)+
\pi\sum_{i=5}^8X_iM_i^{(2)}(\b{X})\]
for appropriate linear forms $M_i^{(j)}$, with
\[H_1(X_5,\ldots,X_8)=X_5X_6+Q_5(X_5,0,X_7,X_8)\]
and
\[H_2(X_5,\ldots,X_8)=X_6X_7+Q_6(X_5,0,X_7,X_8).\]
Thus $G_1$ and $G_2$ will be integral
forms, but since $|\det(U)|_v^2|\det(T)|_v<1$ the pair $G_1,G_2$ will
not be minimized. If we denote the reductions of $H_1$ and $H_2$ 
by $h_1$ and $h_2$ respectively we see that they
are in the right shape for an application of Lemma \ref{7l}, after 
a re-labelling of the variables.  In view of the alternative
conclusions of the lemma there are now two cases to
consider.

\subsection{Case 1.}
If $h_2$ and $h_2$ have
a singular common zero over $F_v$, there will be a change of variable putting
them into the shape
\[h_i(\b{X})=c_iX_5\ell(X_6,X_7,X_8)+w_i(X_6,X_7,X_8)\;\;\;(i=1,2).\]
A further substitution allows us to write
\[h_i(\b{X})=c_iX_5X_6+w_i'(X_6,X_7,X_8)=X_6\lambda_i(X_5,\ldots,X_8)+
w_i'(0,X_7,X_8)\]
for $i=1,2$ and appropriate linear forms $\lambda_1,\lambda_2$.
Hence if we apply the transforms 
\beql{TU2}
T=\diag(1,1,1,1,1,1,\pi,\pi,\pi)\;\;\;\mbox{and}\;\;\;
U=\diag(\pi^{-1},\pi^{-1})
\eeq
the forms $(G_1,G_2)_T^U=(V_1,V_2)$ say, will be integral.  However if
we label the transforms (\ref{TU1}) as $T_1$ and $U_1$, and the transforms
(\ref{TU2}) as $T_2$ and $U_2$, we see that $T_1T_2$ and $U_1U_2$ satisfy
(\ref{TU=}), so that $V_1,V_2$ must be minimized.  On the other hand
one sees that $R(V_1,V_2)\le 6$.  Thus Lemmas \ref{notmin}, \ref{R4}, \ref{R5}
and \ref{R6} show that the pencil $<V_1,V_2>$ contains a form which
splits off three hyperbolic planes, which suffices for the present lemma.

\subsection{Case 2.}

Suppose next that there are at least $5(\# F_v-1)$
non-zero pairs $(a,b)\in F_v^2$ for which $ah_1+bh_2$ is a
sum of two hyperbolic planes. Since $r(q_3,q_4)=3$ or $4$ there are at 
most $4(\# F_v-1)$ non-zero pairs $(a,b)\in F_v^2$ for which
$aq_3+bq_4$ has rank less than 3.  Thus
there is at least one pair such that $ah_1+bh_2$ is a
sum of two hyperbolic planes and $aq_3+bq_4$ splits off a hyperbolic plane.  
Now take $a',b'\in\ov$
to be arbitrary lifts of $a$ and $b$, and consider $G=a'G_1+b'G_2$.
This will take the shape
\[G(\b{X})=\pi G_3(X_1,\ldots,X_4)+G_4(X_5,\ldots,X_8)+
\pi\sum_{i=5}^8X_iM_i(\b{X}),\]
where $\overline{G_3}$ splits off a hyperbolic plane, and
$\overline{G_4}$ is a sum of two hyperbolic planes.

According to Lemma \ref{lift} we may make a change of variable so that
$G$ becomes
\[X_5X_6+X_7X_8+\pi G_5(X_1,\ldots,X_4)\]
with $\overline{G_5}=\overline{G_3}$.  A second application of
Lemma \ref{lift} then shows that $G_5$ splits off a hyperbolic plane, so
that the pencil $<G_1,G_2>$ contains the form $G$ which splits off
three hyperbolic planes.  This establishes Lemma \ref{R7} in Case 2.
\bigskip

We end this section by establishing Lemma \ref{7l}.
We suppose for the proof that $r(q_1,q_2)=4$ since otherwise 
Lemma~\ref{sz} provides a singular common zero. Without loss of
generality we will assume that $q_1(0,0,1)\not=0$.
We may then replace $q_2$ by $q_2+aq_1$ with a suitable value of $a$,
and substitute $Y_2+aY_1$ for
$Y_2$ so as to produce a form in which $q_2(0,0,1)=0$.  We can then write
\[q_2(Y_0,\ldots,Y_3)=Y_2\big(Y_0+\lambda(Y_1,Y_2,Y_3)\big)
+Y_1\ell(Y_1,Y_3).\]
We proceed to substitute $Y_0+\lambda(Y_1,Y_2,Y_3)$ for $Y_0$ so that 
\[q_2(Y_0,\ldots,Y_3)=Y_0Y_2+Y_1\ell(Y_1,Y_3).\]
This transformation puts $q_1$ into the shape
\[q_1(Y_0,\ldots,Y_3)=Y_0Y_1+q_3(Y_1,Y_2,Y_3)\]
for some new quadratic form $q_3$ with $q_3(0,0,1)=q_1(0,0,1)\not=0$.
Now, for any $a\in F$ we have
\[q_1+aq_2=Y_0(Y_1+aY_2)+q_3(Y_1,Y_2,Y_3)+aY_1\ell(Y_1,Y_3)\]
and on substituting $Y_1+aY_2$ for $Y_1$ this becomes
\begin{eqnarray*}
\lefteqn{Y_0Y_1+q_3(Y_1-aY_2,Y_2,Y_3)+a(Y_1-aY_2)\ell(Y_1-aY_2,Y_3)}
\hspace{1cm}\\
&=&\big(Y_0+\mu(Y_1,Y_2,Y_3)\big)Y_1+q_3(-aY_2,Y_2,Y_3)
-a^2Y_2\ell(-aY_2,Y_3).
\end{eqnarray*}
If $\det(q_1+aq_2)\not=0$ this is a sum of two hyperbolic planes provided that 
the binary form $q_3(-aY_2,Y_2,Y_3)-a^2Y_2\ell(-aY_2,Y_3)$ has
a linear factor over $F$.  Since $r(q_1,q_2)=4$ we need to exclude at 
most $4$ values of $a$ for which $\det(q_1+aq_2)$ vanishes.

We claim that, unless the forms $q_1$ and $q_2$ have a singular common
zero over $F_v$, there are at least $9$
values of $a$ such that the polynomial
\[f_a(U)=q_3(-a,1,U)-a^2\ell(-a,U)\]
has a root $u\in F_v$. There will then be at least $5$
values with the additional property that $\det(q_1+aq_2)\not=0$, and then
$bq_1+abq_2$ will be a sum of two hyperbolic planes for any non-zero $b\in
F_v$.  Thus the claim suffices for the proof of Lemma \ref{7l}.

The coefficient of $U^2$ in $f_a(U)$ is $q_3(0,0,1)\not=0$ so that
$f_a(U)$ is quadratic in $U$ for every value of $a$.  Since there are
at most two roots $u$ for any value of $a$ it will
therefore suffice to show that the curve 
\[\cl{X}: f(U,V)=f_V(U)=0\]
has at least 18 affine points over $F_v$.

We first consider the case in which $\ell$ does not vanish
identically, so that $\cl{X}$ is a curve of degree 3.
If $\cl{X}$ is absolutely irreducible the number of projective points over
$F_v$ is at
least $\# F_v+1-2\sqrt{\# F_v}\ge 21$, by the
Hasse--Weil bound.  At most three of these can be at infinity, so that
there are at least $18$ affine points, as required.  If $\cl{X}$ contains a
line defined over $F_v$ there are at least $\#  F_v+1\ge 33$
projective points and we have the same conclusion. There remains the
possibility that $\cl{X}$ splits into three cubic conjugate lines. However
this case cannot arise since $F(U,V)$ contains no term in $U^3$ and a
non-zero term $q_3(0,0,1)U^2$.  This completes the proof in the case
in which $\ell$ does not vanish identically.

In the alternative case in which $\ell$ vanishes identically, 
$q_2$ reduces to $Y_0Y_2$ and $f(U,V)$
becomes $q_3(-V,1,U)$.  If
$\rank(q_3)=3$, the curve $\cl{X}$ has $\# F_v+1$ projective points over
$F_v$, of which at most two lie at infinity.  We may then complete the
argument as before.  Finally, if $\rank(q_3)\le 2$ we can write
\[q_1=Y_0Y_1+q_4\big(\ell_1(Y_1,Y_2,Y_3),\ell_2(Y_1,Y_2,Y_3)\big)\;\;\;
\mbox{and}\;\;\; q_2=Y_0Y_2\]
for suitable linear forms $\ell_1,\ell_2$ over $F_v$.  
If we choose a non-zero point $\b{y}$ over $F_v$, such that
\[y_0=\ell_1(y_1,y_2,y_3)=\ell_2(y_1,y_2,y_3)=0\]
we then see that $\b{y}$ is a singular common zero for $q_1$ and
$q_2$. This completes the proof of the lemma in this second case.

\section{Global Forms Splitting Off 3 Hyperbolic Planes}

Our task now is to deduce Theorem \ref{T1} from Theorem \ref{T2}. In
doing this we will be inspired by the plan outlined by
Colliot-Th\'{e}l\`{e}ne, Sansuc and Swinnerton-Dyer 
\cite[Remark 10.5.3]{CTSSD2}. However our argument looks somewhat
different from that which they proposed.

We begin by producing a global statement related to Theorem \ref{T2}.
\begin{prop}\label{T2g}
Let $k$ be a number field for which every prime ideal has
absolute norm at least 32, and let $Q_1(X_1,\dots,X_8)$ and 
$Q_2(X_1,\dots,X_8)$ be two quadratic forms such that the
projective variety $Q_1(\b{X})=Q_2(\b{X})=0$ is nonsingular and has a
point over every completion $k_v$ of $k$.  Then
there exist $a,b\in k$ such that
$aQ_1+bQ_2$ has rank 8 and splits off at least 3 hyperbolic planes.
\end{prop}

Our first move will be to establish a variant of Theorem \ref{T2} for
archimedean valuations.
\begin{lemma}\label{T2a}
Let $Q_1(X_1,\dots,X_8)$ and 
$Q_2(X_1,\dots,X_8)$ be quadratic forms over $\R$, such that the
projective variety $Q_1(\b{X})=Q_2(\b{X})=0$ is nonsingular.  Then
there is a non-trivial linear combination $aQ_1+bQ_2$ which splits off
at least 3 hyperbolic planes.
\end{lemma}

We prove this by adapting an argument of Swinnerton-Dyer \cite[\S
4]{SD11}.  Let
\[F_{\theta}(z)=\det\big((\sin\theta)Q_1+(\cos\theta)Q_2-zI\big).\]
This is a polynomial with real roots. For any open interval
$I\subseteq\R$ let $n_I(\theta)$ denote the number of roots of
$F_{\theta}(z)$, counted according to multiplicity, lying in
$I$. If $I=(a,b)$ is a finite interval, and $F_{\theta}(z)$ is
non-zero at $z=a$ and $z=b$, we will have
\[n_{(a,b)}(\theta)=
\frac{1}{2\pi i}\int_{\Gamma}\frac{F'_{\theta}(z)}{F_{\theta}(z)}dz,\]
where $\Gamma$ is the path from $a-i$ to $b-i$ to $b+i$ to $a+i$ and
back to $a-i$. This formula makes it clear that there is a
neighbourhood of $\theta$ on which $n_{(a,b)}(\phi)$ is constant,
whenever $F_{\theta}(a)$ and $F_{\theta}(b)$ are non-zero.

We now write $n_+(\theta),n_-(\theta),n_0(\theta)$ for the number of
roots of $F_{\theta}(z)$ which are positive, negative, or zero,
respectively.  It follows from the above that $n_+$ and $n_-$ are
locally non-decreasing.  Moreover, since the variety
$Q_1=Q_2=0$ is nonsingular Lemma \ref{LBP} shows that  
$n_0(\theta)=0$ or 1. We also observe that 
$n_+(\theta)=n_-(\theta+\pi)$ and
$n_-(\theta)=n_+(\theta+\pi)$, so that either $n_+(\theta)$ or 
$n_+(\theta+\pi)$ must be at least 4.  Suppose that $n_+(\theta)\ge 4$
and let
\[\theta_0=\sup\{\xi\in [\theta,\theta+\pi]: \,  n_+(\xi)\ge 3\}.\]
If $n_+(\theta_0)\ge 5$ one must have 
$\theta_0=\theta+\pi$, since $n_+$ is locally non-decreasing to the
right of $\theta_0$.  This
however is impossible since 
\[n_+(\theta+\pi)=n_-(\theta)\le 8-n_+(\theta)\le 4.\]
On the other hand, if $n_+(\theta_0)\le 4$ we have
$n_-(\theta_0)\ge 3$.  Thus there is an interval
$(\theta_0-\delta,\theta_0]$ on which $n_-\ge 3$.   However it 
follows from the definition of $\theta_0$ as a supremum
that there is a point $\phi\in(\theta_0-\delta,\theta_0]$ such that
$n_+(\phi)\ge 3$, whence $\min\big(n_+(\phi),n_-(\phi)\big)\ge 3$.
Thus $(\sin\phi)Q_1+(\cos\phi)Q_2$ has at least 3 positive eigenvalues, and
at least 3 negative ones, so that it will split off at least 3
hyperbolic planes over $\R$. This completes the proof of Lemma \ref{T2a}.
\bigskip

Moving now to our treatment of Proposition \ref{T2g}, we begin by replacing
$Q_1$ and $Q_2$ by suitable scalar multiples, so that they are defined
over $\ok$.  We then define a set $B$ consisting of
all infinite places, all places above 2, and all places
corresponding to prime ideals dividing $\Delta(Q_1,Q_2)$, as given by
(\ref{Fdef}) and (\ref{Ddef}).  

We now claim that for a place $v$ not
belonging to $B$, every nontrivial linear combination $Q=aQ_1+bQ_2$
has $\rank(\overline{Q})\ge 7$, where
$\overline{Q}$ is the reduction to $F_v$, as usual. Then $Q$ must
split off three hyperbolic planes over $k_v$, by Lemma \ref{lift}.

To prove this latter claim we suppose for a contradiction that 
$\rank(\overline{Q})\le 6$ for $Q=aQ_1+bQ_2$
with $a$, say, a unit in $\ov$. After a
suitable change of variable in $\GL_8(\ov)$ we can arrange that
$\overline{Q}$ is a function of $X_1,\ldots,X_6$ only.  Then
$\det(x\overline{Q}+y\overline{Q_2})$ will be divisible by $y^2$.
Thus if $F(x,y)$ is given by (\ref{Fdef}) and
if $\pi$ is a uniformizing element for $k_v$ then $F(x,y)$ will have a
repeated factor modulo $\pi$, contradicting our assumption that
$\pi\nmid\Delta(Q_1,Q_2)$. This suffices for the proof of our claim.
\bigskip

We can now prove Proposition \ref{T2g}.  Since every prime ideal of $\ok$
has absolute norm at least 32 we have $\# F_v\ge 32$ for every finite place
$v$ of $k$.  Then, according to Theorem \ref{T2} and
Lemma \ref{T2a}, for every $v\in B$ there is a non-trivial pair 
$a_v,b_v$ of elements of $k_v$ such that $a_vQ_1+b_vQ_2$ splits off 3
hyperbolic planes in $k_v$. If $v$ is a finite place and $Q=a_vQ_1+b_vQ_2$
then there is a change of variables in $\GL_8(k_v)$ transforming $Q$ 
into $X_1X_2+X_3X_4+X_5X_6+Q'(X_7,X_8)$ with $Q'$ defined over $\ov$.
According to Lemma \ref{lift}, any other quadratic form congruent to
this modulo $\pi$ also splits off 3 hyperbolic planes over $k_v$.  
We therefore deduce
that there is a positive real $\ep_v$ such that $aQ_1+bQ_2$ splits off
3 hyperbolic planes in $k_v$
whenever $|a-a_v|<\ep_v$ and $|b-b_v|<\ep_v$. We may
obtain the analogous conclusion for infinite places by using the local
non-decreasing property for the functions $n_+(\theta),n_-(\theta)$
introduced above.

By weak approximation, there are suitable elements $a,b\in k$
satisfying the additional condition that $\det(aQ_1+bQ_2)\not=0$.
Then to complete the proof we merely
note that a quadratic form splits off three hyperbolic planes over $k$ if and
only if it does so over every completion $k_v$.

\section{The Hasse Principle in the Absence of Small Prime Ideals}

Our next goal is the following major result, which establishes the Hasse
principle for fields which have no prime ideals with norm less than
32.

\begin{prop}\label{T1p}
Let $k$ be a number field for which every prime ideal has
absolute norm at least 32, and let $Q_1(X_1,\dots,X_8)$ and 
$Q_2(X_1,\dots,X_8)$ be two quadratic forms such that the
projective variety 
\[\cl{V}: Q_1(\b{X})=Q_2(\b{X})=0\]
is nonsingular and has a non-trivial point over every completion 
$k_v$ of $k$.  Then there is a non-trivial point over $k$.
\end{prop}

We take $Q$ to be the form $aQ_1+bQ_2$ given by Proposition \ref{T2g}.
Assuming that $a\not=0$, as we may by symmetry, the forms $Q$ and
$Q_2$ then generate the pencil $<Q_1,Q_2>_k$.  The proof of 
Proposition \ref{T1p} now depends on the following lemma.
\begin{lemma}\label{2H}
There is a change of variable in $\GL_8(k)$ such that
\[Q(\b{X})=Q_3(X_1,\ldots,X_6)+X_7L_1(\b{X})+X_8L_2(\b{X})\]
and
\[Q_2(\b{X})=Q_4(X_1,\ldots,X_6)+X_7L_3(\b{X})+X_8L_4(\b{X})\]
with $\rank(Q_3)=4$ and $\rank(Q_4+\alpha Q_3)\ge 5$ for
every $\alpha\in \overline{k}$.  Moreover $Q_3$ and $Q_4$ have a 
nonsingular common zero over every completion $k_v$ of $k$.
\end{lemma}
We will prove this in the next section, but first we show how it
suffices for Proposition \ref{T1p}.

We make a change of variable so that $Q_3$ depends only on
$X_1,\ldots,X_4$.  Then if $Q_5(X_5,X_6)=Q_4(0,\ldots,0,X_5,X_6)$ has 
rank less than 2, or is a hyperbolic plane, it will have a non-trivial 
zero $x_5,x_6$ over $k$, whence $\cl{V}$ will have a non-trivial point at
$(0,\ldots,0,x_5,x_6,0,0)$. We may
therefore assume that $Q_5$ is anisotropic over $k$. We then have
$Q_5(x_5,x_6)=0$ for a point $(x_5,x_6)$ defined only over a quadratic
extension of $k$, so that the variety $Q_3=Q_4$ has a pair of
conjugate singular points.  We now apply the following result,
which is a consequence of Theorem 9.6 of Colliot-Th\'{e}l\`{e}ne, Sansuc and
Swinnerton-Dyer \cite{CTSSD2}.

\begin{lemma}\label{ext}
Let $Q_3(X_1,\ldots,X_6)$ and $Q_4(X_1,\ldots,X_6)$ be quadratic forms
over a number field $k$, such that the projective variety $\cl{Y}: Q_3=Q_4=0$
is absolutely irreducible of codimension 2, and is not a cone.
Suppose that $\cl{Y}$ has a pair of conjugate singular points over $k$, and
assume further that the pencil $<Q_3,Q_4>_{\overline{k}}$ does not
contain two independent forms of rank 4.  Then if $\cl{Y}$ has 
nonsingular points over every completion of $k$ it will have a point over $k$.
\end{lemma}

We proceed to investigate the possibility that $\cl{Y}$ might be a cone, 
or might fail to be absolutely irreducible of codimension 2.
It is trivial that $\cl{Y}$ has a point over $k$ 
if $\cl{Y}$ is a cone.  Moreover if
$\cl{Y}$ is not absolutely irreducible of codimension 2 then, according to
Lemma \ref{arc}, either $r_{\min}(Q_3,Q_4)\le 2$ or  $r(Q_3,Q_4)\le 4$, in
the notation of \S \ref{PQF}.  In the first case
we would have $a_1Q_3+a_2Q_4=R(X_1,X_2)$ say, and hence
\[a_1Q+a_2Q_2=R(X_1,X_2)+X_7L'(X_1,\ldots,X_8)+X_8L''(X_1,\ldots,X_8),\]
say. This however has rank at most 6 in contradiction to Lemma
\ref{LBP}. Thus we cannot have
$r_{\min}(Q_3,Q_4)\le 2$.  On the other hand, if
$r(Q_3,Q_4)\le 4<6$ then Lemma \ref{lemshape1} provides a
$k$-point, given by taking $x_1=\ldots=x_r=0$ and $x_R=1$
in the notation of the lemma.

Thus in any cases in which Lemma \ref{ext} is not applicable we will
automatically have a point on $\cl{Y}$ defined over $k$.  Thus Proposition
\ref{T1p} follows from Lemma \ref{ext}.

\section{Proof of Lemma \ref{2H}}

Lemma \ref{2H} requires us to control both local solvability
and ranks, and we begin by investigating the local solvability condition.
Our first result is the following.

\begin{lemma}\label{ls}
Let $S_1(X_1,\ldots,X_n)$ and $S_2(X_1,\ldots,X_n)$ be quadratic forms
over $k$ such that $r_{\min}(S_1,S_2)\ge 5$ and $r(S_1,S_2)\ge 7$.  
Then there is a finite
set $B=B(S_1,S_2)$ of places of $k$ such that for any $v\not\in
B$, and any linear form $L(X_1,\ldots,X_n)$ defined over $k$,
the variety 
\[\cl{V}_L:\, S_1=S_2=L=0\]
has a nonsingular point over $k_v$.
\end{lemma}

Without loss of generality we may assume that $S_1,S_2$ and $L$ are
defined over $\ok$. Moreover, if $L(X_1,\ldots,X_n)=\sum \ell_i X_i$
we can multiply $L$ by a suitable constant so that the ideal generated
by $\ell_1,\ldots, \ell_n$ has norm at most $c_k$, the Minkowski
constant. We shall require $B$ to include all the infinite places,
together with all finite places corresponding to prime ideals of norm
at most $c_k$.

When $v\not\in B$ at least one coefficient, $\ell_n$ say, of $L$ is a unit
in $\ov$.  We may then replace the variety $\cl{V}_L$ by
$S_{1,L}(X_1,\ldots,X_{n-1})=S_{2,L}(X_1,\ldots,X_{n-1})=0$, where
\begin{eqnarray*}
\lefteqn{S_{i,L}(X_1,\ldots,X_{n-1})}\hspace{2cm}\\
&=&S_i\big(X_1,\ldots,X_{n-1},
-\ell_n^{-1}(\ell_1X_1+\ldots,\ell_{n-1}X_{n-1})\big), \;\;\;(i=1,2).
\end{eqnarray*}
and ask whether the variety $S_{1,L}=S_{2,L}=0$ has a nonsingular
point in $k_v$.

It follows from the Lang--Weil theorem \cite{LW}
that there is a constant $C$ depending only on $n$, such that the variety
\[\overline{\cl{V}_L}:\,\overline{S_{1,L}}=\overline{S_{2,L}}=0\]
has a nonsingular point
over $F_v$ provided firstly that it is absolutely irreducible
of codimension 2, and secondly that $\# F_v\ge C$.  By Hensel's Lemma
these conditions then suffice for the existence of a nonsingular point on
$\cl{V}_L$ over $k_v$. We therefore require $B$ to include all places 
corresponding to prime ideals of norm at most $C$, and all places such
that there is some linear form $L$ for which the variety
$\overline{\cl{V}_L}$ fails to be absolutely 
irreducible of codimension 2 over $F_v$.
In view of Lemma \ref{arc} it is sufficient to add places for 
which there is some linear from $L$ having either
\beql{lsc1}
r_{\min}(\overline{S_{1,L}},\overline{S_{2,L}})\le 2
\eeq
or
\beql{lsc2}
r(\overline{S_{1,L}},\overline{S_{2,L}})\le 4.
\eeq

Thus $\cl{V}_L$ will always have a nonsingular point over $k_v$ when
$v\not\in B$, and it remains to prove that $B$ is finite.
We begin by considering those $v$ for which (\ref{lsc1}) holds for 
some $L$.  Thus there exist $a,b\in
F_v$, not both zero, and forms $q,\ell_1,\ell_2$ over $F_v$, such that
\[a_1\overline{S_{1,L}}+b\overline{S_{2,L}}=\ell_1\ell_2.\]
We then deduce that there is a third form $\ell_3$ such that
\[a_1\overline{S_{1}}+b\overline{S_{2}}=\ell_1\ell_2+\overline{L}\ell_3.\]
Thus the $5\times 5$ minors of the matrix
$M(X\overline{S_{1}}+Y\overline{S_{2}})$, which are binary forms in
$X$ and $Y$, all vanish at $X=a, Y=b$.  However the $5\times 5$ minors of
$M(XS_{1}+YS_{2})$ cannot have a common zero, since
$r_{\min}(S_1,S_2)\ge 5$.  Thus there are two of them which have a
non-zero resolvent ${\rm Res}\in\ov$ say.  This gives us an element depending
only on $S_1$ and $S_2$. However the minors can only have a common
zero over $F_v$ when $\overline{{\rm Res}}=0$, and this can only happen for a
finite set of places $v$.

In the alternative case for which (\ref{lsc2}) holds, any linear
combination $a_1\overline{S_{1}}+b\overline{S_{2}}$ differs from the
corresponding form $a_1\overline{S_{1,L}}+b\overline{S_{2,L}}$ by a
multiple of $\overline{L}$, so that
\[\rank\big(a_1\overline{S_{1}}+b\overline{S_{2}}\big)\le
\rank\big(a_1\overline{S_{1,L}}+b\overline{S_{2,L}}\big)+2\le 4+2\]
for every pair $a_1,a_2\in\overline{F_v}$.  Hence all $7\times 7$
minors of $M(X\overline{S_{1}}+Y\overline{S_{2}})$ must vanish identically.
However the $7\times 7$ minors of
$M(XS_{1}+YS_{2})$ cannot all vanish identically, since
$r(S_1,S_2)\ge 7$.  Thus there is some non-zero coefficient $\mu$ of
some minor, such that (\ref{lsc2}) can hold only when
$\overline{\mu}=0$.  This too can hold only for a finite set of places
$v$. We have therefore shown that we can take $B$ to be finite,
thereby completing the proof of Lemma \ref{ls}.
\bigskip

In our next result we handle the remaining places.
\begin{lemma}\label{ls2}
Let $S_1(X_1,\ldots,X_n)$ and $S_2(X_1,\ldots,X_n)$ be quadratic forms
over $k$ such that the variety $S_1=S_2=0$ has a nonsingular point over
every completion $k_v$. Suppose that
$r_{\min}(S_1,S_2)\ge 5$ and $r(S_1,S_2)\ge 7$, and that
$S_1$ splits off two hyperbolic planes. 

If $P$ is a nonsingular projective point on $S_1=0$, defined over $k$, let
$L_P(X_1,\ldots,X_n)=0$ be the tangent hyperplane to $S_1=0$ at $P$.
Then there is a Zariski-dense set of such $P$ for which the variety
\[\cl{V}(P): S_1(X_1,\ldots,X_n)=S_2(X_1,\ldots,X_n)=L_P(X_1,\ldots,X_n)=0\]
has a nonsingular point over every completion $k_v$.
\end{lemma}

Although the lemma requires $S_1=S_2=L_P=0$ to have nonsingular
solutions for every place $v$ it is clear that Lemma \ref{ls} may be
applied, producing a finite set $B$ of places outside which a
nonsingular solution is guaranteed for any linear form $L_P$.

For each $v\in B$ we have a nonsingular point
$\b{x}_v=(x_{1v},\ldots,x_{nv})\in k_v^n$ on $S_1=S_2=0$. Since
$\b{x}_v$ is  nonsingular we have $M(S_1)\b{x}_v\not=\b{0}$.  We
now consider the codimension 2 quadric defined over $k_v$ by 
\[\cl{Q}: S_1(\b{X})=\b{X}^tM(S_1)\b{x}_v=0.\]
Since $S_1$ splits off two hyperbolic planes over
$k$ we may write it as 
\[S_1(\b{X})=X_1X_3+X_2X_4+S_3(X_5,\ldots,X_n),\]
after a change of variable. We then choose non-zero points
\[\b{y}=(a,b,0,\ldots,0)\;\;\;\mbox{and}\;\;\;\b{z}=(0,0,c,d,0,\ldots,0)\]
over $k_v$,
both lying on the hyperplane $\b{X}^tM(S_1)\b{x}_v=0$.  Thus both
$\b{y}$ and $\b{z}$ lie on the quadric $\cl{Q}$.  Moreover 
\[\nabla S_1(\b{y})=(0,0,a,b,\ldots,0)\;\;\;\mbox{and}\;\;\;
\nabla S_1(\b{z})=(c,d,0,\ldots,0).\]
Since these cannot both be
proportional to $M(S_1)\b{x}_v$ we deduce that $\cl{Q}$ has at least one 
nonsingular point, $\b{w}_v$ say, over $k_v$. When an irreducible quadric has a
nonsingular point over an infinite field such points are automatically
Zariski-dense.  Thus in our case the available points $\b{w}_v$ 
cannot be restricted to a line.
We may therefore suppose that our point $\b{w}_v$ is chosen so that
$M(S_1)\b{w}_v$ does not lie on the line through $M(S_1)\b{x}_v$ and 
$M(S_2)\b{x}_v$.

We therefore have
\[S_1(\b{x}_v)=S_2(\b{x}_v)=S_1(\b{w}_v)=0,\]
\[\b{w}_v^tM(S_1)\b{x}_v=0,\]
\[\nabla S_1(\b{w}_v)\not=\b{0},\]
and
\beql{lsmat}
\rank\left(\begin{array}{ccc}
\frac{\partial S_1(\b{x}_v)}{\partial X_1} & \ldots &
\frac{\partial S_1(\b{x}_v)}{\partial X_n} \\
\rule{0mm}{5mm}\frac{\partial S_2(\b{x}_v)}{\partial X_1} & \ldots &
\frac{\partial S_2(\b{x}_v)}{\partial X_n} \\
\rule{0mm}{5mm}\frac{\partial S_1(\b{w}_v)}{\partial X_1} & \ldots &
\frac{\partial S_1(\b{w}_v)}{\partial X_n} \\
\end{array}\right)=3.
\eeq
Thus $P_v=\b{w}_v$ is a nonsingular point on $S_1=0$, and $\b{x}_v$ is
a nonsingular point on the variety $\cl{V}(P_v)$ defined in the lemma.

For each place $v$ there is a $3\times 3$ determinant 
$\Delta(\b{x}_v,\b{w}_v)$
formed from the matrix (\ref{lsmat}) such that
$\Delta(\b{x}_v,\b{w}_v)\not=0$.  If $|\Delta(\b{x}_v,\b{w}_v)|_v
=\delta_v$ then there
is a positive real $\ep_v$ such that $|\Delta(\b{x}_v,\b{w})|_v=\delta_v$
whenever $|\b{w}-\b{w}_v|<\ep_v$. We now use weak approximation on the
quadric $Q_1=0$ to choose a point $P=\b{w}$ suitably close to $\b{w}_v$
for each $v\in B$.  Since $\b{w}_v$ is a nonsingular point on $Q_1=0$
the resulting point $P$ will also be nonsingular. We then claim that,
if $\b{w}$ is sufficiently close to $\b{w}_v$, the variety $\cl{V}(P)$ will
have a nonsingular point over $k_v$.  This will follow from Hensel's
Lemma, using the starting value $\b{x}=\b{x}_v$, for which 
$|\Delta(\b{x}_v,\b{w})|_v=\delta_v>0$ if $\b{w}$ is close enough to
$\b{w}_v$. Moreover we have $Q_1(\b{x}_v)=Q_2(\b{x}_v)=0$, while
$L_P(\b{x}_v)$ can be made sufficiently small for Hensel's Lemma to
apply, merely by taking $\b{w}$ suitably close to $\b{w}_v$. Thus we
do indeed obtain a point on $\cl{V}(P)$, and since the lifting argument for
Hensel's Lemma preserves the value of $|\Delta(\b{x},\b{w})|_v$ we see
that the resulting point is nonsingular.  This establishes Lemma
\ref{ls2}, since we can use any point $\b{w}$ sufficiently close to
$\b{w}_v$ for each $v\in B$.
\bigskip

We are now ready to prove Lemma \ref{2H}, which will require us to use
Lemma \ref{ls2} twice.  
The hypotheses of Lemma \ref{ls2} are
satisfied for the forms $S_1=Q$ and $S_2=Q_2$, since
$r_{\min}(Q,Q_2)\ge 7$. Let us write temporarily $\cl{I}$ for the variety
$S_1=0$, and $\cl{W}$ for the variety $S_1=S_2=0$, both considered in
$\proj^7$. Each of these varieties is nonsingular. We now proceed to consider
their dual varieties $\cl{I}^*$ and $\cl{W}^*$. The
reader should recall that the dual of a variety $\cl{Y}$ is the closure 
of the set of hyperplanes which are tangent at a nonsingular 
point of $\cl{Y}$. For any variety $\cl{Y}\subset \proj^m$, 
the dual $\cl{Y}^*$ is a proper subvariety of $\proj^{m*}$, and is
irreducible.   Moreover we have $(\cl{Y}^*)^*=\cl{Y}$.
In our case $\cl{I}^*$ will be a quadric hypersurface.  We claim
that $\cl{I}^*$ cannot be contained in $\cl{W}^*$.  Indeed since
$\cl{I}^*$ is an irreducible hypersurface and $\cl{W}^*$ is an irreducible
proper subvariety of $\proj^7$ the only situation in which one could have
$\cl{I}^*\subseteq \cl{W}^*$ is when $\cl{I}^*=\cl{W}^*$.
However this would imply that
\[\cl{I}=(\cl{I}^*)^*=(\cl{W}^*)^*=\cl{W}.\]
The claim then follows since $\cl{I}\neq\cl{W}$.  

It now follows from Lemma \ref{ls2} that we may choose our point $P$ so
that the hyperplane $L_P=0$ is not tangent to $\cl{W}$.  Thus the variety
$\cl{V}(P)$ will be nonsingular. Since $P$ is defined over $k$ the
hyperplane $L_P=0$ is also defined over $k$. We may therefore make a
change of variables in $\GL_8(k)$ so that the hyperplane $L_P=0$ is $X_8=0$. 
The forms $S_1$ and $S_2$ then become
\[S_i(\b{X})=T_i(X_1,\ldots,X_7)+X_8L_i(X_1,\ldots,X_8),\;\;\;(i=1,2)\]
say, with $\rank(T_1)=6$ by Lemma \ref{add}. The conclusions of
Lemma \ref{ls2} imply that the forms $T_1$ and $T_2$ 
have a nonsingular common solution over every completion $k_v$. 
Moreover the variety $T_1=T_2=0$ will be nonsingular
since $\cl{V}(P)$ is nonsingular. In particular Lemma \ref{LBP} implies that 
$r_{\min}(T_1,T_2)=6$ and $r(T_1,T_2)=7$.
Moreover, as $S_1=Q$ splits off three hyperbolic planes we 
deduce that $T_1$ splits off two hyperbolic planes.  

We therefore see that the forms $T_1$ and $T_2$ satisfy the hypotheses
for a second application of Lemma \ref{ls2}. Since $r(T_1,T_2)=7$ there
are at most 7 linear combinations $W=T_2+\alpha T_1$ with
$\alpha\in\kb$ such that $\rank(W)<7$. Let these be $W_1,\ldots,W_m$,
say.  For each of these the variety
$\cl{X}_i: W_i=0$ is distinct from $\cl{X}: T_1=0$ so that, 
by the argument above, their
duals, which are irreducible quadric hypersurfaces, are different. 

We may therefore choose $P$ so that the hyperplane
$L_P=0$ belongs to $\cl{X}^*$ but to none of the sets $\cl{X}_i^*$.  We now
repeat the manoeuvres above.  We make a change of variables in
$\GL_7(k)$ so that the hyperplane $L_P=0$ is $X_7=0$.  The forms
$T_1,T_2$ then become
\[T_i(X_1,\ldots,X_7)=U_i(X_1,\ldots,X_6)+X_7M_i(X_1,\ldots,X_7),\;\;\;(i=1,2)\]
say.  Lemma \ref{add} shows that $\rank(U_1)=\rank(T_1)-2=4$, while
Lemma \ref{ls2} tells us that $U_1$ and $U_2$
have a nonsingular common solution over every completion $k_v$. 

We claim now that $\rank(U_2+\alpha U_1)\ge 5$ for every $\alpha\in
\kb$.  Since our change of variables has transformed $(T_2+\alpha
T_1)(X_1,\ldots,X_7)$ into
\beql{UT}
(U_2+\alpha U_1)(X_1,\ldots,X_6)+X_7L_{\alpha}(X_1,\ldots,X_7)
\eeq
for a suitable linear form $L_{\alpha}$, we see that
$\rank(U_2+\alpha U_1)\ge\rank(T_2+\alpha T_1)-2$ for every $\alpha$.
This verifies the claim except when $T_2+\alpha T_1$ is one of the
forms $W_i$ above.  Suppose then that $T_2+\alpha T_1=W_i$.  Since
$r_{\min}(T_1,T_2)\ge 6$ we know that $\rank(W_i)=6$.
However, if $\rank(U_2+\alpha U_1)\le 4$ then there are independent
linear forms $\lambda_1(X_1,\ldots,X_6),\ldots,
\lambda_6(X_1,\ldots,X_6)$ say, such that one can
write $U_2+\alpha U_1$ as a form in $\lambda_1,\ldots,\lambda_4$
alone. Thus (\ref{UT}) produces
\[W_i=U(\lambda_1,\ldots,\lambda_4)+X_7L(\lambda_1,\ldots,\lambda_6, X_7),\]
say.  Since $\rank(W_i)=6$ the linear form $L(\lambda_1,\ldots,\lambda_6, X_7)$
must properly contain at least one of $\lambda_5$ or $\lambda_6$.
One then sees from (\ref{UT}) that the variety $\cl{X}_i: W_i=0$ has
a nonsingular point with $\lambda_1=\ldots=\lambda_4=X_7=0$ and $L\not=0$.
The tangent hyperplane at this point would be $X_7=0$, which is also the
tangent hyperplane of $\cl{X}: T_1=0$ at $P$. This however is impossible, since
$P$ was chosen so that the hyperplane $L_P=0$ was in none of the dual
varieties $\cl{X}_i^*$.  This completes the proof of our claim.

We now see that the effect of our two applications of Lemma \ref{ls2}
is to produce a change of variables putting $Q$ and $Q_2$ into the
shape
\[Q(\b{X})=U_1(X_1,\ldots,X_6)+X_7L_1(\b{X})+X_8L_2(\b{X})\]
and
\[Q_2(\b{X})=U_2(X_1,\ldots,X_6)+X_7L_3(\b{X})+X_8L_4(\b{X})\]
for suitable linear forms $L_1(\b{X}),\ldots,L_4(\b{X})$, all defined 
over $k$.  Moreover we have arranged that the variety $U_1=U_2=0$ has
nonsingular points in every completion $k_v$ of $k$, that
$\rank(U_1)=4$ and that $\rank(U_2+\alpha U_1)\ge 5$ for every $\alpha\in
\kb$.  This therefore suffices for Lemma \ref{2H}, and thereby also
completes our demonstration of Proposition \ref{T1p}.

\section{Deduction of Theorem \ref{T1}}

In this final section we first show how to remove the condition that all
prime ideals of $\ok$ have norm at least 32, and then explain why the
weak approximation property holds automatically in our situation.

We begin by choosing a prime $q\ge 5$ which does not divide 
$[k:\Q]$, and we proceed
to construct a number field $\Q(\theta)$ of degree $q$ over $\Q$,
such that every prime ideal of $\mathcal{O}_{\Q(\theta)}$ has norm at least
32. Clearly it suffice that every rational prime $p\le 31$ is inert in 
$\Q(\theta)/\Q$.
For each such prime we choose a monic polynomial
$f_p(X)\in\Z[X]$ of degree $q$ which is irreducible modulo $p$.  We then use the
Chinese Remainder Theorem to produce a monic polynomial $f(X)\in\Z[X]$
 of degree $q$, with $f(X)\equiv f_p(X)\modd{p}$ for every $p\le 31$.  
Then $f$ is
irreducible over $\Q$ since it is irreducible modulo 2, for example.
We claim that $\Q(\theta)$ will be a suitable field, where $\theta$ is a root of
$f(X)$.  Let $p\le 31$ be prime.  Then $f_p$ is irreducible modulo 
$p$, since it has no repeated factors over $\mathbb{F}_p$.  Thus
$p\nmid{\rm Disc}(\theta)$, so that we may apply the Kummer--Dedekind 
theorem to deduce that $p$ is inert in $\Q(\theta)$. It follows that $N(P)\ge
  32$ for every prime ideal of $\mathcal{O}_{\Q(\theta)}$, as
  required.  Indeed we have $N_{\Q(\theta)/\Q}(I)\ge 32$ for every
  non-trivial integral ideal $I$.

We now consider the field $k'=k(\theta)$. If $P$ is a prime ideal of
  $\mathcal{O}_{k'}$ then
\[N_{k'/\Q}(P)=N_{\Q(\theta)/\Q}\big(N_{k'/\Q(\theta)}(P)\big)\ge 32.\]
We also note that $[k':\Q]=[k':\Q(\theta)][\Q(\theta):\Q]$, and
$[k':\Q]=[k':k][k:\Q]$, so that both $[\Q(\theta):\Q]=q$ and
$[k:\Q]$ divide $[k':\Q]$. It follows that $q[k:\Q]$ divides
$[k':\Q]$, since we chose $q$ to be coprime to $[k:\Q]$.  On the other
hand it is clear that $[k(\theta):k]\le q$, whence $[k':\Q]\le
q[k:\Q]$.  We therefore conclude that $[k':\Q]=q[k:\Q]$, so that
$[k':k]=q$.

Since $\cl{V}$ has points over every completion of $k$ it will also
have points over every completion of $k'$.  Thus Proposition \ref{T1p} is
applicable, and shows that $\cl{V}$ has a point $\x=(x_1,\ldots,x_8)$
say, over $k'$.  It follows that the quadratic form
$Q_1(\b{X})+TQ_2(\b{X})$, which is defined over the function field
$k'(T)$, has a non-trivial point at $\x$.  However
$Q_1(\b{X})+TQ_2(\b{X})$ is also defined over the subfield $k(T)$, and
$[k'(T):k(T)]=[k':k]=q$, which is odd.  It therefore follows from a
  result of Springer \cite{springer} that $Q_1(\b{X})+TQ_2(\b{X})$ has
  a zero over $k(T)$. Finally, we apply the Amer--Brumer Theorem 
\cite[\S 3, Satz 7]{amer}, \cite{brumer}, which we state as follows.
\begin{lemma}\label{AB}
Let $S_1(X_1,\ldots,X_n)$ and $S_2(X_1,\ldots,X_n)$ be quadratic forms
defined over a field $K$ of characteristic not equal to 2. Then if
$S_1+TS_2$ has a non-trivial zero over the function field $K(T)$, the
forms $S_1$ and $S_2$ have a simultaneous zero over $K$.
\end{lemma}
We therefore conclude that our forms $Q_1$ and $Q_2$ have a
simultaneous zero over $k$.
\bigskip

It remains to show that $\cl{V}$ has the weak approximation property.
This however is a direct consequence of Theorem 3.11 of Colliot-Th\'el\`ene,
Sansuc and Swinnerton-Dyer \cite{CTSSD1}.

\bigskip
\bigskip

Mathematical Institute,

24--29, St. Giles',

Oxford

OX1 3LB

UK
\bigskip

{\tt rhb@maths.ox.ac.uk}

\end{document}